\documentclass{hkart}
\usepackage{amsmath,amsfonts,amssymb,color,enumerate}

\def\cen{\centerline}
\def\nin{\noindent}
\def\defi{:=}
\def\hat{\widehat}

\def\<{\langle} \def\>{\rangle}
\def\[{[[} \def\]{]]}
\def\.{{\cdot}}
\def\u#1{\underline#1}
\def\ssk{\smallskip}
\def\til{\widetilde}
\def\phi{\varphi}
\def\SI{{\cal S}}
\def\cA{{\cal A}} 
\def\cB{{\cal B}} 
\def\dA{{\cal A}_d} 
\def\dB{{\cal B}_d} 
\def\SI{{\cal S}} 
 
\def\V{{\cal V}}
\def\W{{\cal W}}
\def\AA{{\cal W}\hskip-2pt {\cal A}}
\def\WAC{{\cal W}\hskip-2pt {\cal A}\hskip-.5pt{\cal C}}

\def\LL{{\cal W}\hskip-1pt {\cal L}}
\def\TA{{\cal T}\hskip-3pt {\cal A}}
\def\II{{\cal I}}
\def\Ii{{\mathbb I}}

\def\JJ{{\cal J}}
\def\li #1{#1_{\rm Lie}}
\def\un #1{#1_{\rm u}}
\def\fin #1{#1_{\rm fin}}
\def\WA{\AA}
\def\fWA{\AA_{\rm fin}}
\def\fW{\W_{\rm fin}}

\def\C{\mathbb C}
\def\c{\mathfrak c}

\def\i{\mathfrak i}

\def\s{\mathfrak s}

\def\G{\mathbb G} 
\def\I{\mathbb I}
\def\I{{\bf I}}
\def\N{\mathbb N}
\def\P{\mathbb P} 
\def\PP{\mathfrak P} 
\def\Q{\mathbb Q}
\def\R{\mathbb R}
\def\S{\mathbb S}

\def\compgr{\mathcal CO\hskip-2pt M\hskip-2pt PG\hskip-1ptR}
\def\progr{\mathcal P\hskip-2pt ROG\hskip-1ptR}
\def\prolie{\mathcal P\hskip-2pt RO\hskip-1ptLI\hskip-1pt E}
\def\K{\mathbb K} 
\def\L{\mathfrak L}
\def\g{\mathfrak g}
\def\h{\mathfrak h}
\def\i{\mathfrak i}
\def\T{\mathbb T}
\def\Z{\mathbb Z}

\def\SS{{\bf S}}

\def\Span{\mathop{\rm span}\nolimits}

\def\incl{\mathop{\rm incl}\nolimits}

\def\Hom{\mathop{\rm Hom}\nolimits}

\def\id{\mathop{\rm id}\nolimits}

\def\im{\mathop{\rm im}\nolimits}
\def\ob{\mathop{\rm ob}\nolimits}

\def\TT{\mathop{\bf T\hphantom{}}\nolimits}
\def\UU{\mathop{\bf U\hphantom{}}\nolimits}
\def\VV{\mathop{\bf V\hphantom{}}\nolimits}
\def\bsk{\bigskip} \def\msk{\medskip}
\def\lead{\leaders\hbox to 1.5ex{\hss${.}$\hss}\hfill}
\def\arr{\hbox to 40pt{\rightarrowfill}}
\def\larr{\hbox to 40pt{\leftarrowfill}}
\def\mapdown#1{\Big\downarrow\rlap{$\vcenter{\hbox{$\scriptstyle#1$}}$}}
\def\lmapdown#1{\llap{$\vcenter{\hbox{$\scriptstyle#1$}}$}\Big\downarrow}
\def\mapright#1{\smash{\mathop{\arr}\limits^{#1}}}
\def\lmapright#1{\smash{\mathop{\arr}\limits_{#1}}}
\def\mapleft#1{\smash{\mathop{\larr}\limits^{#1}}}
\def\lmapleft#1{\smash{\mathop{\larr}\limits_{#1}}}
\def\mapup#1{\Big\uparrow\rlap{$\vcenter{\hbox{$\scriptstyle#1$}}$}}
\def\lmapup#1{\llap{$\vcenter{\hbox{$\scriptstyle#1$}}$}\Big\uparrow}

\def\a{{\bf a}} \def\b{{\bf b}} \def\1{{\bf1}}
\def\aa{{\mathfrak a}} \def\bb{{\mathfrak b}}
\newcommand{\End}{\mathrm{End}}

\newcommand{\eps}{\varepsilon}

\title{On Weakly Complete Universal Enveloping Algebras:\\ 
        A Poincar\'e--Birkhoff--Witt Theorem}
                  
\author{Karl Heinrich Hofmann and Linus Kramer\thanks{
Both authors were supported by Mathematisches Forschungsinstitut Oberwolfach
in the program RiP (Research in Pairs). Linus Kramer is
funded by the Deutsche Forschungsgemeinschaft
under Germany's
Excellence Strategy EXC 2044-390685587,
Mathematics M\"unster: Dynamics-Geometry-Structure.}}

\lastname{Hofmann and Kramer}



\address{%
Karl Heinrich Hofmann\\       
Fachbereich Mathematik\\      
Technische Universit\"at Darmstadt\\
Schlossgartenstra{\ss}e 7\\
64289 Darmstadt, Germany\\       
\hbox{hofmann@mathematik.tu-darmstadt.de}}

\address{%
Linus Kramer\\	
Mathematisches Institut\\ 
Universit\"at  M\"unster\\
Einsteinstra{\ss}e 62\\
48149 M\"unster, Germany\\
\hbox{linus.kramer@uni-muenster.de}}

\volume{32}
\annum{2022}
\editor{Karl-Hermann Neeb}
\received{December 3, 2021} 
\final{January 17, 2022} 
\begin{document}
\maketitle
\vglue-10pt

\begin{abstract}
The {\sc Poincar\'e-Birkhoff-Witt} Theorem deals with the  
structure and universal property of the universal enveloping
algebra $U(L)$ of a Lie algebra $L$, e.g., over $\R$ or $\C$.  
In $\<1\>$``K.H.~Hofmann and L.~Kramer, {\it On Weakly Complete Group Algebras 
of Compact Groups}, J. Lie Theory {\bf30} (2020), 407--426'',  
the weakly complete universal enveloping algebra $\UU(\g)$ of a profinite 
dimensional topological Lie algebra $\g$ was introduced. Here
it is shown that the classical universal enveloping algebra $U(|\g|)$
of the  abstract Lie algebra underlying $\g$ is a dense subalgebra
of $\UU(\g)$, algebraically generated by $\g\subseteq \UU(\g)$.
It is further shown that, inspite of $\UU$ being  a left 
adjoint functor, it nevertheless  
preserves projective limits in the form 
$\UU(\lim_\i \g/\i)\cong \lim_\i\UU(\g/\i)$,
for profinite-dimensional Lie algebras $\g$ represented as 
projective limits of their finite-dimensional quotients. The required
theory  is presented in an appendix which is of
independent interest.---In a natural
way, a weakly complete enveloping algebra $\UU(\g)$ is 
 a weakly complete symmetric Hopf algebra with a 
Lie subalgebra $\P(\UU(\g))$ of {\it primitive} elements containing $\g$
(indeed properly if $\g\ne\{0\}$), and with a nontrivial multiplicative
pro-Lie group $\G(\UU(\g))$ of {\it grouplike} units,
having $\P(\UU(\g))$ as its Lie algebra--in  contrast with the 
classical  {\sc Poincar\'e-Birhoff-Witt} environment of $U(L)$, 
thus providing a new aspect of Lie's Third Fundamental
Theorem: Indeed a canonical pro-Lie subgroup $\Gamma^*(\g)$ of $\G(\UU(\g))$
is identified whose Lie algebra is naturally isomorphic to $\g$.
The structure of $\UU(\g)$ is  described in detail for $\dim\g=1$. 
The  primitive and grouplike components and their mutual relationship  
are evaluated precisely.---In $\<1\>$ and 
$\<2\>$``R. Dahmen and K.H.~Hofmann, 
The Pro-Lie Group Aspect of Weakly Complete Algebras$(\dots)$,
J. of Lie Theory {\bf29} (2019), 413--455'', 
the real weakly complete group Hopf algebra
$\R[G]$ of a compact group $G$ was described. In particular, the set
$\P(\R[G]))$ of primitive elements of $\R[G]$ was identified 
as the Lie algebra $\g$ of $G$.  It is now shown that for 
any compact group $G$ with Lie algebra $\g$ there  is a natural
morphism of weakly complete symmetric Hopf algebras  
$\omega_\g\colon\UU(\g)\to\R[G]$, implementing the identity on $\g$ and
inducing a morphism of pro-Lie groups $\Gamma^*(G)\to\G(\R[G])\cong G$:
yet another aspect of Sophus Lie's Third Fundamental Theorem!
\end{abstract}

\section{The Weakly Complete Enveloping Algebra\\ 
of a Profinite-Dimensional Lie Algebra}

In \cite{hofkra} we have initiated the theory of 
{\it  weakly complete universal enveloping algebras over $\K$}
hoping that  in some fashion  this concept would resemble the 
classical universal enveloping
algebra of a Lie algebra such as it is presented in the famous
{\sc Poincar\'e-Birkhoff-Witt}-Theorem 
(see e.g.\ \cite{bour}, Chap. 1,  Paragraph 2, n$^{\rm o}$ 7,
Th\'eor\`eme 1., p.30). While this was not exactly the case,
we shall discuss now how close we come to that theorem.

So we let $\K$  denote one of the topological fields $\R$ or $\C$.
For a topological Lie algebra $\g$ over $\K$ we let $\II(\g)$ denote the filter
basis of all closed ideals $\i\subseteq\g$ such that $\dim \g/\i<\infty$.

\begin{Definition} \label{1.1a} A topological Lie algebra $\g$ over $\K$ is 
called {\it profinite-dimensional} if $\g=\lim_{\i\in\II(\g)}\g/\i$.
Let $\LL$ denote the category of profinite-dimensional Lie algebras
(over $\K$) and continuous Lie algebra morphisms between them.
\end{Definition}

Notice that by its definition every profinite-dimensional 
Lie algebra is weakly complete.
A comment following 
Theorem 3.12 of \cite{dhtwo}  exhibits an example of a weakly complete
$\K$-Lie algebra which is not a profinite-dimensional Lie algebra.

Let $\AA$ denote the category of weakly complete associative 
unital algebras over
$\K$. However, instead of considering the full category 
of weakly complete Lie algebras over $\K$, in the following 
we consider $\LL$, the category of profinite-dimensional Lie algebras over
$\K$ and continuous $\K$-Lie algebra morphisms.
The reason for this restriction is Theorem \ref{1.1}
stating that every weakly complete unital $\K$-algebra is the
projective limit of its finite-dimensional quotient algebras. 
This implies at once the following 
\begin{Proposition} \label{1.2} Let $A$ be any weakly complete 
unital $\K$-algebra and $\li A$ 
the weakly complete Lie algebra obtained
by considering on the weakly complete vector space $A$ the Lie
algebra obtained with the Lie bracket $[x,y]=xy-yx$. Then $\li A$
is profinite-dimensional.
\end{Proposition} 

The functor  which associates with a weakly 
complete associative algebra $A$ the profinite-dimensional 
Lie algebra $\li A$  is
called the {\it underlying Lie algebra functor}. 

\medskip
For the complete proof of the following existence theorem,
we shall invoke a considerable portion of a bulk  category theoretical
arguments.  We shall collect these in an appendix, since,
firstly, they reach far beyond the current application and,
secondly, their full presentation might have led the reader astray from
the present line of thought had we presented them at this point.

\begin{Theorem}  \label{1.3a} {\rm (The Existence Theorem of $\UU$)} 
The underlying Lie algebra functor $A\mapsto\li A$
from $\WA$ to $\LL$ has a left adjoint $\UU\colon\LL\to\AA$.

The front adjunction $\lambda_\g\colon\g\to\li{\UU(\g)}$
is an embedding of profinite-dimension\-al Lie algebras. 
\end{Theorem}

\begin{Proof} The category $\LL$ is complete. 
(Exercise. Cf.\ Theorem A3.48 of 
\cite{compbook}, p.\ 819.) The
 ``Solution Set Condition'' (of Definition A3.58 in \cite{compbook}, 
p.\ 824) holds.
(Exercise: Cf.\ the proof Lemma 3.58 of \cite{compbook}, p.\ 91.) 
Hence $\UU$ exists by the Adjoint Functor Existence Theorem
(i.e., Theorem A3.60 of \cite{compbook}, p.\ 825).

The assertion about $\lambda_\g$ being an embedding follows from 
Proposition \ref{morphs} (ii) in the Appendix.
\end{Proof}

In other words, 
each profinite-dimensional Lie algebra $\g$ may be considered
as a closed Lie subalgebra of
$\li{\UU(\g)}$ with the property  that  each continuous Lie algebra
 morphism $f\colon \g\to \li A$ for some weakly complete associative unital 
algebra $A$ extends uniquely to a $\WA$-morphism $f'\colon \UU(\g)\to A$. 

$$
\begin{matrix}& \LL&&\hbox to 7mm{} &\AA\cr 
\noalign{\vskip3pt}
\noalign{\hrule}\cr
\noalign{\vskip3pt}%
   \g&\mapright{\lambda_\g}&\li{\UU(\g)}&\hbox to 7mm{} &\UU(\g)\\
\lmapdown{\forall f}&&\mapdown{\li{(f')}}&\hbox to 7mm{}&
         \mapdown{\exists! f'}\\
 \li A&\lmapright{\id}&\li A&\hbox to 7mm{}&A.
\end{matrix}
$$

\medskip

\noindent If necessary we shall write $\UU_\K$ instead of $\UU$ 
whenever the ground field should be emphasized.

\begin{Definition} \label{1.4a} For each profinite-dimensional $\K$-Lie
algebra, we shall call $\UU_\K(\g)$ 
{\it the weakly complete enveloping algebra} of $\g$ 
(over $\K$).
\end{Definition}

\medskip
\begin{Remark} \label{3.X}  For every profinite-dimensional Lie 
algebra $\g$, a morphism $f\colon\g\to\K$, $f(x)=0$, 
according to the definition of $\UU_\K(\g)$ induces
a natural $\WA$-morphism $\alpha_\g\colon \UU_\K(\g)\to\K$ such that 
$\alpha_\g(\g)=\{0\}$ and that $\alpha_\g\circ\iota_{{\UU}(\g)}=\id_\K$.     
\end{Remark}
\msk

The retraction  $\alpha_\g$ is also called the {\it augmentation}
of $\UU_\K(\g)$. 

\msk
In the Appendix we shall also introduce for any weakly complete 
vector space $W$ its weakly complete tensor algebra $\TT(W)$
(cf.\ paragraph ({\bf C}) preceding Proposition \ref{morphs} and
Theorem \ref{1.3}) and show that $\UU(\g)$ is a quotient algebra
of $\TT(|\g|)$ if $|\g|$ is the underlying weakly complete vector space
underlying $\g$.
There is a commutative diagram
$$ \begin{matrix}\TT(|\g|)&\mapright{\alpha_{|\g|}}&\K\\
            \lmapdown{f'} && \mapdown{=}\\
              \UU_\K(\g)&\lmapright{\alpha_\g}&\K.\end{matrix}$$

\msk
Moreover, we have the following corollary to our existence theorem:

\medskip
\begin{Corollary} \label{density3} 
For any profinite-dimensional Lie algebra $\g$,
the unital associative subalgebra $\<\g\>$ generated algebraically
in $\UU(\g)$ by $\g$ is dense in $\UU(\g)$. 
\end{Corollary}

\begin{Proof} The assertion follows from Proposition \ref{epimorphs}
in the Appendix.
\end{Proof}

\bsk
Of course we would like to have a better  insight into the structure
of the algebra $\<\g\>$. This information we provide in the following
section and thereby close the gap between the concepts of the
weakly complete enveloping algebra and the classical universal enveloping
abstract algebra dealt with in the Poincar\'e-Birkhoff-Witt Theorem.

\section{The Abstract Enveloping Algebra \mbox{$U(L)$} 
of a Lie Algebra \mbox{$L$}}

We briefly recall that the functor which assigns to a unital $\K$-algebra
$X$ the underlying Lie algebra $\li X$ (with the underlying $\K$
vector space of $X$ as vector space structure  endowed with
the bracket operation $(x,y)\mapsto [x,y]\defi xy-yx$ as Lie bracket)
has a left adjoint functor $U$ which assigns to a Lie algebra $L$
a unital associative algebra $U(L)$ and  a natural
Lie algebra morphism $\rho_L\colon L\to \li{U(L)}$ such that
for each Lie algebra morphim $f\colon L\to \li X$ for a unital
algebra $X$ there is a unique morphism of unital algebras 
$f'\colon U(L)\to X$ such that $f=\li {f'}\circ \rho_L$.
The algebra $U(L)$ is called the 
{\it universal enveloping algebra of}  $L$.
A large body of text book literature is available on it.
A prominent result is the Poincar\'e-Birkhoff-Witt Theorem
on the structure of $U(L)$ which implies in particular
that $\rho_L\colon L\to\li{U(L)}$ is injective.

\bsk

From the Theorem of Poincar\'e, Birkhoff and Witt it is known
that $\rho_L$ is injective. One may therefore assume that 
$L\subseteq U(L)$ such that $\rho_L$ is the inclusion function. 
(See \cite{bour} or \cite{dix}.)
In this parlance the universal property reads as follows:

\ssk

\noindent
{\it For each unital algebra $A$,  each Lie algebra morphism
$f\colon L\to \li A$ extends uniquely to an algebra
morphism $f'\colon U(L)\to A$.}

\ssk
 Also from the
Theorem of Poincar\'e, Birkhoff and Witt we know that 

\ssk

\noindent
{\it $U(L)$ is the unital algebra generated by $L$, i.e.,
$U(L)=\<L\>$.}

\bsk
In the present section we shall now denote by $|\g|$ the abstract
Lie algebra underlying the profinite-dimensional Lie algebra
$\g$. Then the main result of this section will be a 
complete clarification of  the relation of the weakly complete
enveloping algebra $\UU(\g)$ of a
profinite-dimensional Lie algebra $\g$ and the universal
enveloping algebra $U(|\g|)$ of  $|\g|$.

\medskip

\begin{Lemma} \label{eps-1} For a profinite-dimensional Lie algebra $\g$
there is a natural morphism $\eps_\g\colon U(|\g|)\to |\UU(\g)|$
of unital algebras such that
\begin{enumerate}[{\rm(i)}]
\item the following diagram is commutative:
$$\begin{matrix} |\g|&\mapright{\incl}&U(|\g|)\\
\lmapdown{\id_{|\g|}}&&\mapdown{\eps_\g}\\
|\g|&\lmapright{|\incl|}& |\UU(\g)|.\end{matrix}$$ 
\item The image of $\eps_\g$ is dense in $\UU(\g)$.
\item The morphism $\eps_\g$ is injective if $\g$ is finite-dimensional.
\end{enumerate}
\end{Lemma}
\begin{Proof} (i) The claim is a direct consequence of the universal
property of the functor $U$.

(ii) We have $\im(\eps_\g)=\eps_\g(U(|\g|))=\eps_\g(\<|\g|\>
=\<\eps_\g(\g)\>=\<|\g|\>$ in $|\UU(\g)|$. From Corollary \ref{density3}
we know that $\<|\g|\>$ is dense in $\UU(\g)$.

(iii) If $\g$ is finite-dimensional, then every finite-dimensional 
Lie algebra representation  $\rho\colon \g\to \li A = \li{\End(V)}$
for a finite-dimensional vector space $V$ extends to an associative
representation $\rho'\colon \UU(\g)\to \End(V)$ . Then 
$\rho\circ \eps_\g\colon U(\g)\to\End(V)$ is an extension to an
associative representation of $U(\g)$ which is unique.
By Harish-Chandra's Lemma (see Dixmier \cite{dix}, 2.5.7),
the extensions of associative representations
of $U(\g)$ of finite-dimensional Lie algebra  representations of $\g$
separate the points of $U(\g)$ and so the claim follows.
\end{Proof}

The remainder of this section now is devoted to removing the restriction to
finite-dimensionality in Lemma \ref{eps-1}(iii). That is, we want to
show
 
\begin{Lemma} \label{eps-2} For a profinite-dimensional Lie algebra $\g$,
the algebra morphism $\eps_\g\colon U(|\g|)\to |\UU(\g)|$ is injective. 
\end{Lemma}

The proof will occupy the remainder of this section.
We shall resort to
the existing literature on $U(L)$ such as \cite{bour} or \cite{dix}.
We are given the profinite-dimensional Lie algebra $\g$ and we 
write $L\defi|\g|$ for the underlying Lie algebra. 
So $L\subseteq \li{U(L)}$. 
Let $B$ be a totally ordered basis of $L$.
 We begin by recalling the following basic fact
from the Poincar\'e-Birkhoff-Witt Theorem (see \cite{bour},
Corollary 3, Section 7 of Paragraph 2):
$$\til B{\defi}\{b_1b_2{\cdots}b_m| 1{\le}m,\ b_1,\dots,b_m\in B,\ 
b_1{\le}b_2{\le}\cdots{\le}b_m\}\mbox{ is  a basis of $U(L)$.}
\leqno{\rm(PBW)}$$
Now assume that 
$$0\ne u\in U(L).$$
Then there is a finite subset $F\subseteq B$ such that
$u\in \Span(\til F)$ for 
$$\til F=\{b_1b_2\cdots b_m| 1\le m,\quad b_1,\dots,b_m\in  F,\quad
         b_1\le b_2\le\cdots\le b_m\}\subseteq\til B.$$

\begin{Lemma} \label{einsA} There is a closed ideal $J$ of $\g$
so that $J\cap \Span F= \{0\}$.
\end{Lemma}
\begin{Proof}
The vector space $V=\Span F$ is finite-dimensional.  
Let $C$ be the boundary of a compact 0-neighborhood  in $V$. Then
$V=\K\.C$. 

\ssk

Returning at this point to the fact that $\g$ is 
a profinite-dimensional Lie algebra, we conclude that 
there is a filterbasis $\II$
of  closed ideals $I$ of $\g$ such that $\dim\g/I<\infty$ 
for $I\in\II$ such
that $\g\cong\lim_{I\in\II}\g/I$. In particular, 
$\bigcap \II=\{0\}$. Therefore
$\bigcap_{I\in\II}(C\cap I)=C\cap\{0\}=\emptyset$. Since $C$ is compact,
the filter basis $\{C\cap I:I\in\II\}$ with empty intersection
consists of compact sets  and therefore must contain the empty set. 
Thus there is a $J\in \II$ such that $C\cap J=\emptyset$. 
In fact, we have $J\cap\Span F=\{0\}$. Indeed,
suppose there were 
a nonzero $t\in\K$ and a
$c\in C$ such that $t\.c\in J$, then $c=t^{-1}\.(t\.c)\in t^{-1}\.J=J$
which is  impossible.
\end{Proof} 

\ssk

\begin{Lemma} \label{eins} Assume that there is  
an ideal $J$ of $L$ such that  $J\cap\Span F=\{0\}$.
Then the image of $u$ under
the morphism $U(L)\to U(L/J)$ is nonzero.
\end{Lemma}

\begin{Proof} 
We choose a finite dimensional vector subspace $H$ of $L$
containing $F$ such that $L=H\oplus J$. Let $E$ be a totally ordered
basis for $H$ such that the order of $E$ extends that of $F$, choose
a totally ordered basis $D$ of $J$ and make sure that $E\cup D$
has a total order extending the orders of $E$  and $D$, thus
yielding a totally ordered basis  of $L$.

We consider the quotient morphism of Lie algebras $q_J\colon L\to L/J$.
Then $q_J$ maps $E$ bijectively onto a basis $E'=q_J(E)$ of $L/J$.
Let $$\til E=\{b_1\cdots b_m\, |\quad m\ge1, b_1,\dots,b_m\in E\}.$$ Then 
$U(q_J)\colon U(L)\to U(L/J)$ maps $\til E$ bijectively onto a basis
$\til{E'}$ of $U(L/J)$  by (PBW). If we now write 
$u=\sum_{S\in \til E}c_S\.S$ with $c_S\in \K$, then
$U(q_J)(u)=\sum_{q_J(S)\in\til{E'}} c_S\.q_J(S)\ne0$, 
since $\til{E'}$ is a basis of $U(L/J)$. 
\end{Proof}

As the kernel of the quotient map $U(L)\to U(L/J)$ is $U(L)J$,
the claim of the preceding lemma may be expressed equivalently
in the form $u\notin U(L)J$.

Now we recall that $q_J\colon L\to L/J\subseteq U(L/J)$ 
is in fact the underlying Lie algebra
morphism $|p_J|$ of a quotient morphism $p_J\colon\g\to \g/J$
of profinite Lie algebras and that $q_J$
 extends  uniquely 
to an algebra morphism $U(|p_J|)\colon U(L)\to U(L/J)$
with kernel $U(L)J$. Then from Lemma \ref{eins} we know that
$U(|p_J|)(u)$ is nonzero in $U(L/J)$ and from Lemma \ref{eps-1}
we infer that $\eps_{\g/J}$ is injective. The commutative diagram
$$\begin{matrix}|\g|&\mapright{\incl}&U(|\g|)&\mapright{\eps_\g}&|\UU(\g)|\\ 
\lmapdown{|p_J|}&&\lmapdown{U(|p_J|)}&&\mapdown{|\UU(p_J)|}&&\\
|\g/J|&\mapright{\incl}&U(|\g/J|)&\mapright{\eps_{\g/J}}&|\UU(\g/J)|\end{matrix}$$
then shows that $\eps_\g(u)\ne0$. Therefore, since $u\in U(L)\setminus\{0\}$
was arbitrary, $\eps_\g$ is injective, leaving the elements of $L=|\g|$ fixed.
This completes the proof of Lemma \ref{eps-2}. Thus
$U(|\g|)$ may be considered as a subalgebra of $|\UU(\g)|$, containing 
$|\g|\subseteq|\UU(\g)|$.

This may be rephrased in the following Theorem which summarizes our efforts
to elucidate the close relation between $U(|\g|)$ and   $\UU(\g)$:

\begin{Theorem} \label{main-1} {\rm (The Relation of $U(-)$ and $\UU(-)$)} 
For any profinite-dimensional real or complex Lie algebra
$\g$ considered as a closed Lie subalgebra of $\li{\UU(\g)}$, the associative unital 
subalgebra $\<\g\>$ generated algebraically by $\g$ in $\UU(\g)$ is naturally isomorphic
to $U(|\g|)$ (under an isomorphism fixing the elements of $\g$) and is dense
in $\UU(\g)$.
\end{Theorem}

In a slightly careless sense we may memorize this as saying:

\nin
{\it For a profinite-dimensional Lie algebra $\g$, 
the weakly complete topological
enveloping algebra $\UU(\g)$ is ``a completion of $U(|\g|)$'', and  we have
$$\g\subseteq \<\g\>=U(\g)\subseteq \overline{U(\g)}=\UU(\g).\leqno(*)$$}

\section{The Projective Limit Preservation\\of the Weakly Complete 
Enveloping  Functor {\bf U}}

Since every weakly complete unital algebra is a strict projective
limit of all finite-dimensional quotient algebras, it will now turn out to
be sufficient to test the universal property of the functor $\UU$ only for
{\it finite-dimensional} unital associative algebras:

\begin{Proposition} \label{finitedim2} Assume that the profinite-dimensional Lie algebra
$\g$ is contained functorially in a weakly complete unital algebra $\VV(\g)$
such that for each \emph{finite-dimensional} unital algebra $A$ and
each morphism of profinite-dimensional Lie algebras $f\colon\g\to\li A$
there is a unique morphism of weakly complete unital algebras
$f'\colon \VV(\g)\to A$ extending $f$. Then $\VV(\g)\cong \UU(\g)$
naturally.
\end{Proposition}

\begin{Proof} We apply  the Density and 
Adjunction Theorem \ref{density} in the Appendix with $\cA$ as the
category of weakly complete associative unital algebras, and $\cB$
as the category of profinite-dimensional Lie algebras with the full
subcategory $\cB_d$ of finite-dimensional Lie algebras which is
topologically dense in $\cB$. Then, by
hypothesis,
the function $\VV\colon\ob\cB\to\ob\cA$ is conditionally left
adjoint to the functor $(\cdot)_{\rm Lie}\colon \cA\to \cB$ which
maps an associative algebra to the Lie algebra with the Lie
bracket $[a,b] =ab-ba$ (see Definition \ref{basic}). 
Then by  Theorem \ref{density}, $\VV$
is naturally isomorphic to the left adjoint $\UU$ 
of $(\cdot)_{\rm Lie}$.
\end{Proof}

Perhaps more deeply we shall see now that, while
as a left-adjoint functor, $\UU$ preserves colimits,
is also preserve certain limits, namely, the projective
limits $\g=\lim_{\i\in\II(\g)}\g/\i$ of Definition \ref{1.1a}.
Indeed
in Theorem \ref{solution-b} in the Appendix we  show:

\begin{Theorem} \label{lim2} {\rm($\UU$ preserves some projective limits)} 
For a profinite-dimensional Lie algebra $\g$ with its filter basis $\II(\g)$
of cofinite-dimensional ideals $\i$ we have
$$\g\cong\lim_{\i\in\II(\g)}\g/\i\mbox{ in }\LL \mbox{ and } 
\UU(\g)\cong\lim_{\i\in\II(\g)}\UU(\g/\i)\mbox{ in }\AA.$$
\end{Theorem}

The argument in the Appendix shows, that while the assertion of
the theorem is natural and easy to absorb, its proof is deeper
than one would expect initially.

\section{The Weakly Complete Universal Enveloping Algebra\\ 
         as a Hopf Algebra}

We now address the important aspect of enveloping algebras from their beginning,
namely, the fact that they are symmetric Hopf algebras. 
For some of the proofs in this section  we refer to our predecessor paper \cite{hofkra}.

\begin{Proposition} \label{8.3}  
The universal enveloping functor $\UU$ is multiplicative,
that is, there is a natural isomorphism 
$\UU(\g_1\times \g_2)\to \UU(\g_1)\otimes_\W\UU(\g_2)$.
\end{Proposition}

\nin
For a proof see \cite{hofkra}, Proposition 6.3.

\begin{Lemma} \label{li-morph} For any weakly complete unital algebra
$A$, the vector space morphism $\Delta_A\colon A\to A\otimes_\W A$,
$\Delta_A(a)=a\otimes 1+1\otimes a$
is a morphism of weakly complete 
Lie algebras  $\li A\to \li{(A\otimes_WA)}$.
\end{Lemma}

\nin
Cf. \cite{hofkra}, Lemma 6.4.

\msk

\nin
Recall the natural morphism $\lambda_\g\colon\g\to\li{\UU(\g)}$ which we
consider as an inclusion morphism.  By Lemma \ref{li-morph},

\centerline{%
$p_\g=\delta_{\UU(\g)}\circ \lambda_\g\colon \g\to 
\li{(\UU(\g)\otimes_\W\UU(\g))}$}

\nin is a morphism of
weakly complete  Lie algebras. By 
the universal property of $\UU$, $p_\g$ yields
a unique natural morphism of weakly complete associative unital
algebras
$\gamma_\g\colon \UU(\g)\to \UU(\g)\otimes_\W \UU(\g)$ 
such that $p_\g=\li{(\gamma_\g)}\circ \lambda_\g$. 
Recall the augmentation $\alpha_\g\colon \UU_\K(\g)\to\K$ 
(see Remark \ref{3.X})
and the inclusion morphism $\iota_{\UU(\g)}\colon \K\to \UU_\K(\g)$
(see Remark \ref{Remark-2.X}). Accordingly, we have an idempotent endomorphism
$$ \iota_{\UU_K(\g)}\circ \alpha_\g\colon \UU_\K(\g)\to \UU_\K(\g).$$
Further, the augmentation $\alpha_\g$ acts as coidentity, and the function
$x\mapsto -x:\UU(\g)\to \UU(\g)$ as symmetry as is readily checked for
$x\in \g$, and $\g$ generates $\UU_\K(\g)$ as topological algebra
by Corollary \ref{density3}.

Now we have

\begin{Proposition} \label{hopf} {\rm ($\UU(\g)$ as a Hopf algebra)} 

\noindent{\rm(a)} Each weakly complete enveloping algebra $\UU(\g)$ is
a weakly complete symmetric Hopf algebra with the comultiplication 
$\gamma_\g$ and the augmentation  $\alpha_\g\colon\UU(\g)\to\K$ as coidentity.

\noindent{\rm(b)} If $f\colon\g\to\h$ is a morphism of profinite-dimensional
Lie algebras, then the morphism $\UU_\K(f)\colon\UU_K(\g)\to\UU_K(\h)$ respects
comultiplication, coidentity, and symmetry, that is, $\UU_\K(f)$ is a morphism
of symmetric Hopf algebras. 
\end{Proposition}

\begin{Proof} For (a), see \cite{hofkra}, Corollary 6.5.

\msk

For (b) we consider a morphism $f\colon\g\to\h$ a morphism of 
profinite-dimensional Lie algebras and for the functoriality 
of $\UU$ (short for $\UU_K(-)$) as regards to comultiplication
$\gamma_\g\colon \UU(\g)\to \UU(\g)\otimes\UU(\g)$
we verify the commutativity of the following diagram (with $\otimes=\otimes_\W$)
$$\begin{matrix}
\UU(\g)&\mapright{\gamma_\g}&\UU(\g)\otimes\UU(\g)&\mapright{\cong}&\UU(\g\times\g)\\
\lmapdown{\UU(f)}&&\lmapdown{\UU(f)\otimes\UU(f)}&&\mapdown{\UU(f\times f)}\\
\UU(\h)&\lmapright{\gamma_\h}&\UU(\h)\otimes\UU(\h)&\lmapright{\cong}&\UU(\h\times\h).
\end{matrix}$$
(See also Proposition \ref{8.3}.)
Coidentity and symmetry are treated similarly.
\end{Proof}

 \medskip 

This proposition expresses the fact that $\UU_\K$ is a functor 
from the category of profinite-dimensional Lie algebras to 
the category of weakly complete symmetric Hopf algebras.                          
Its significance is emphasised by
 the fact that  essential portions of the noteworthy theory
of weakly complete symmetric Hopf algebras have meanwhile entered
the textbook literature. 
(See \cite{compbook}, Appendix A3, Appendix A7, Chapter 3--Part 3.)
We have collected some essential features in our Appendix such as
Theorems \ref{expone} and \ref{aug}. 

Now we specialize these to the case of $A=\UU_\K(\g)$. 
We use the notation
$\R_<=\{r\in\R:0<r\}$ and recall that $\g\subseteq A$
and that $A^\times$ denotes the group of units of $A$. 
For the exponential function 
$\exp\colon\li A\to A^\times$  as in Theorem \ref{expone} 
we define the closed subgroup 
$$\Gamma^*(\g)\defi\overline{\<\exp \g\>}\subseteq A^\times.$$ 
The following theorem now is a principal result in the theory
of weakly complete enveloping algebras of profinite-dimensional
real or complex Lie algebras.

\begin{Theorem} \label{8.5} {\rm(The Weakly Complete Enveloping Hopf Algebra)}
 Let $\g$ be a profinite-dimensional Lie algebra and $\UU(\g)$ its  
weakly complete  enveloping algebra containing $\g$ according to
{\rm Theorem \ref{main-1}}. 
Then the following statements hold:
\begin{enumerate}[\rm(a)]

\item The group of units $\UU(\g)^\times$ is dense in $\UU(\g)$. 
It is an almost connected  pro-Lie group,
connected in the case of $\K=\C$.
The algebra $\UU(\g)$ has an exponential function 
$\exp\colon{\li{\UU(\g)}}\to \UU(\g)^\times$. 
The Lie algebra $\L(\UU(\g)^\times)$  of $\UU(\g)^\times$ is
(naturally isomorphic to) $\li{\UU(\g)}$.

\item  The pro-Lie algebra $\P(\UU(\g))$ 
is the Lie algebra of the pro-Lie group
$\G(\UU(\g))$  of grouplike elements and 
the restriction and corestriction of 
$\exp$ is the exponential function for
this group.

\item The profinite-dimensional Lie algebra $\P(\UU(\g))$ contains 
$$\g=\P(U(|\g|))=\P(\UU(\g))\cap U(|\g|).$$ 
 
\item For $\K=\R$, the restriction and corestriction of $\exp$ yields
the exponential function 
$$\exp_{\Gamma^*(\g)}\colon \g=\L(\Gamma^*(\g))\to\Gamma^*(\g) $$ 
of that pro-Lie subgroup $\Gamma^*(\g)$ of $\UU_\R(\g)^\times$ whose
Lie algebra is precisely $\g$.

\item Define the hyperplane  ideal $\Ii$ as the kernel of the augmentation 
$\alpha_\g$. Then we have
\begin{enumerate}[\rm(i)]
\item for $\K=\R$:\ $\exp(\UU_\R(\g))=(\R_<\.1)\oplus \Ii$, an open half space, 
\item for $\K=\C$:\ $\exp(\UU_\C(\g))=\UU_\C(\g)\setminus\Ii$.
\end{enumerate}
\end{enumerate}
\end{Theorem}

\begin{Proof} For the proofs of (a) and (b)  see Theorem \ref{expone}
in the Appendix. 

The proof of  (c) follows from
\cite{hofkra}, Theorem 3.4 and Theorem \ref{main-1}. Cf.\ also \cite{compbook},
Theorem A3.102 and its proof for $\K\in\{\R, \C\}$. 

The proof of (d)
must verify that $\L(\overline{\<\g\>}) \subseteq \g$. This conclusion we derive 
from \cite{probook}, Corollary 4.22 and its proof.  

The proof of (e) follows from Theorem \ref{aug} in the Appendix. 
\end{Proof}

Here are some immediate consequences:

\begin{Corollary}  \label{8.6} 
The weakly complete enveloping algebra $\UU(\g)$ of a
nonzero  pro\-finite-dimensional weakly complete 
Lie algebra $\g$  has nontrivial grouplike elements
contained in  $\UU(\g)^\times$.
Specifically, there is  a pro-Lie subgroup $\Gamma^*(\g)$
of grouplike elements
whose Lie algebra is isomorphic to $\g$ and whose exponential
function is induced by that of $\UU(\g)$. 
\end{Corollary} 

By contrast, on the purely algebraic side, the universal enveloping Hopf 
algebra $U(L)$ of a Lie
algebra $L$ shows no visible nontrivial grouplike elements
while a nontrivial weakly complete enveloping algebra always does. 

\nin
We shall see that even in the case of the smallest possible nonzero
candidate $\g=\K$, the space $\P(\UU(\g))$ is substantially larger
than $\g$ (see Theorem \ref{monoth} below).
In the classical setting of the discrete enveloping Hopf algebra
in characteristic 0 we have $\P(U(L))=L$: 
see e.g.\ \cite{serre-1}, Theorem 5.4
on p.~LA 3.10. 

\msk

\begin{Corollary} \label{Lie-3-1} For any profinite-dimensional
Lie algebra $\g$ there is a pro-Lie group $G$ whose Lie algebra
$\L(G)$ may be identified with $\g$.
\end{Corollary}

Indeed the theorem provides a weakly complete unital algebra $A$ 
with an exponential function $\exp_A$  such that 
$\exp_G\colon \L(G)\to G$ may be identified with a restriction
and corestriction of $\exp_A$.
 
This is indeed much more than  what is historically known
as {\it Sophus Lie's Third Fundamental Theorem}.
\bsk

\subsection{Lie's Third Fundamental Theorem for 
profinite-dimensional Lie algebras}

It is worthwhile to elucidate the insight that our present context
throws a new light on Lie's Third Fundamental Theorem. Therefore we
recall the contemporary aspect of this background:
\ssk

\begin{Theorem} \label{lie3} {\rm(Sophus Lie's Third Principal Theorem)}  
For every profinite-di\-mens\-ional
real Lie algebra $\g$ there is a simply connected pro-Lie group 
$\Gamma(\g)$, whose Lie algebra $\L(\Gamma(\g))$ is (isomorphic to) 
$\g$. For any pro-Lie group $G$
with Lie algebra $\g$ there is a quotient morphism 
$\alpha_\g\colon \Gamma(\g)\to G$
such that the following diagram commutes:
$$\begin{matrix}\g&\mapright{=}&\g\\
\lmapdown{\exp_{\Gamma(\g)}}&&\mapdown{\exp_G}\\
\Gamma(\g)&\lmapright{\alpha_\g}&G.\end{matrix}$$
\end{Theorem}

For a systematic proof see  \cite{liethree}, 
or e.g.\ \cite{probook}, Chapter 6, p.~249, see notably Theorem 6.4, p.~232. 
Our  Theorem \ref{lie3}
is also cited in \cite{compbook}, Theorem A7.29. For the definition of simple
connectivity see \cite{compbook}, Definition A2.6. 
Let us recall here that for an abelian $\g$ (that is, a weakly complete real 
vector space), the underlying vector space of $\g\cong \L(\Gamma(\g))$ 
is isomorphic to $\Gamma(\g)$ via 
$\exp_{\L(\Gamma(\g))}\colon\L(\Gamma(\g))\to \Gamma(\g)$.

\msk

Theorem \ref{lie3}  applies at once to $G=\Gamma^*(\g)$ as follows:

\begin{Corollary} For each profinite-dimensional real Lie algebra
$\g$ there is a natural morphism $\alpha_\g\colon \Gamma(\g)\to \Gamma^*(\g)$
such that the diagram  $$\begin{matrix}\g&\mapright{=}&\g\\
\lmapdown{\exp_{\Gamma_\g}}&&\mapdown{\exp_{\Gamma^*(\g)}}\\
\Gamma(\g)&\lmapright{\alpha_\g}&\Gamma^*(\g)\end{matrix}$$
is commutative.\end{Corollary}

\msk
The pro-Lie group $A^\times$ of units of
a weakly complete unital associative algebra $A$ has the property
that finite-dimensional continuous representations separate
the points, and so any pro-Lie group injected into such a group
$A^\times$ shares this property. Consider the Lie algebra
$\g={\rm sl}(2,\R))$ of the Lie group $G={\rm SL}(2,\R)$, 
and let $\til G=\Gamma(\g)$ be 
the universal covering group of $G$. Every continuous linear 
representation of $\til G$, however, factorizes through $G$
 (see\ \cite{hilnee}, p.590, Example 16.1.8), and therefore
$\Gamma(\g)$ cannot be injected into any group of the form 
$A^\times$. Hence for $\g={\rm sl}(2,\R)$ the morphism
$\alpha_\g\colon\Gamma(\g)\to\Gamma^*(\g)$ cannot be injective
and thus certainly cannot be an isomorphism.
 
\bsk
 
\section{The Abelian Case} 

For an abelian Lie algebra $\g$, the weakly complete unital
algebra $\UU(\g)$ is commutative. In various special aspects
we considered this situation 
in \cite{dhtwo}, Lemmas 3.4, 3.5, and 3.10ff., and
in \cite{hofkra}, Section 5 and Example 6.2. We now return to 
the commutative situation more systematically now and
discuss the structure of $\UU_\K(\g)$ completely for $\dim \g=1$,
and derive consequences for the abelian case in general.

\subsection{The power series algebra}

\ssk\nin A first and simplest step is the discussion of the power 
series algebra.
We recall the notation $\N=\{1,2,3,\dots\}$ and $\N_0=\{0,1,2,3\dots\}$.
The set $\N_0$ is a semiring for addition and multiplication (that is, an
addition and multiplicative commutative monoid with distributivity).

\begin{Lemma} The weakly complete vector space $\S\defi\K^{\N_0}$ 
supports a monoid multiplication called {\it convolution} as follows:
Let $\a{=}(a_k)_{k\in \N_0}{\in}\S$ and $\b{=}(b_m)_{m\in\N_0}{\in}\S$. Then
$$\a*\b=\Big(\sum_{k+m=n}a_kb_m\Big)_{n\in\N_0}. \leqno(1)$$ 
With pointwise addition and convolution, $\S=(\S,+,*)$ is a 
weakly complete topological algebra.
\end{Lemma}

\begin{Proof} The verification is an exercise:

\ssk\nin We write $X=(0,1,0,0,\dots)$ and observe 
$$\begin{matrix} X&=&(0,1,0,0,0,\dots)\\
                 X^2&=&(0,0,1,0,0,\dots)\\
                 X^3&=&(0,0,0,1,0,\dots)\\
                 \vdots&& \vdots\end{matrix}$$
Then
$$\a*X=(a_n)_{n\in \N}*X=a_0\.X+a_1\.X^2+ a_2\.X^3 +\cdots.$$
The projective limit representation 
$\S{\cong}\lim_{n\in\N}\frac{\K[X]}{X^{n+1}\K[X]}$
completes the proof.
\end{Proof}

Accordingly, $\S$ is called the {\it power series algebra in one variable},
usually written as  $\K[[X]]$. 
 
\ssk

It is useful to  recall that in the category $\W$ of weakly complete $\K$-vector
spaces, for any pair of sets {\bf X} and {\bf Y} we have a natural isomorphism
$$\K^{\bf X}\otimes_\W\K^{\bf Y}\to\K^{{\bf X}\times {\bf Y}}\leqno(*)$$
induced by the bijection 
$(a_x)_{x\in\bf X}\otimes(b_y)_{y\in \bf Y})
\mapsto(a_xb_y)_{(x,y)\in{\bf X}\times{\bf Y}}$.

\ssk

We have a multiplication $\mu\colon \S\otimes_\W\S\to \S$
according to 
$$\mu({\bf a}\otimes {\bf b})=
\Big(\sum_{k+m=n}a_kb_m\Big)_{n\in\N_0}.$$
We write  $\1\defi(1,0,0,\dots)$ and set 
$$X_1\defi X\otimes\1
     =(0,1,0,\dots)\otimes(1,0,0,\dots)\in\S\otimes_\W\S,$$
and 
$$X_2\defi \1\otimes X
     =(1,0,0,\dots)\otimes(0,1,0,\dots)\in\S\otimes_\W\S.$$ 
So we obtain $X_1X_2=X_2X_1$ in $\S\otimes_\W\S$ and compute 
$$\Big(\sum_{m\in\N_0}a_mX_1^m\Big)\Big(\sum_{n\in\N_0}b_nX_2^n\Big)=
\sum_{m,n\in\N_0}a_mb_nX_1^m X_2^n$$
for $\a=(a_m)_{m\in\N_0}$ and 
$\b=(b_n)_{n\in\N_0}$ in $\S$.
Thus 
$$\S\otimes_\W\S=
\Big\{\sum\nolimits_{(m,n)\in\N_0\times\N_0}c_{mn}X_1^mX_2^n:c_{mn}\in\K\Big\}$$
is the ring of power series in two commuting variables. 
Write  $\a=\sum_{n\in\N_0}a_nX^n\in\S$ and note 
${\bf a}\otimes{\bf 1} = \sum_{m\in\N_0}a_mX_1^m$ and 
${\bf 1}\otimes{\bf a} = \sum_{n\in\N_0}a_nX_2^n$  in $\S\otimes_\W\S$.
We then have two morphisms of vector spaces
 $$\Delta, \gamma\colon \S\to \S\otimes_\W\S\mbox{ as follows:}$$
 For $\a=\sum_{n\in\N_0}a_nX^n$ 
$$\Delta(\a)\defi\a\otimes\1+\1\otimes\a                                                 
=\sum_{m,n\in\N_0}(a_mX_1^m+a_nX_2^n)\leqno(2)$$
 and
$$\gamma(\a)\defi\gamma\Big(\sum_{n\in\N_0}a_n\.X^n\Big)
     =\sum_{n\in\N_0}a_n(X_1+X_2)^n,\leqno(3)$$
where $\gamma$ is in fact a morphism of weakly complete
algebras. 
Also, there is an
identity $\epsilon\colon\K\to\S$, $\epsilon(t)=t\.\1.$ and
a coidentity (or {\it augmentation} $\kappa\colon \S\to\K$
given by $\kappa(\a)=\kappa((a_n)_{n\in\N_0})\defi a_0$,
 and a symmetry $\sigma\colon\S\to\S$
given by $\sigma(X)\defi-X$, that is 
$\sigma(\a)=\sigma((a_n)_{n\in\N_0})\defi ((-1)^na_n)_{n\in\N}$.
The following diagram  is commutative:
$$\begin{matrix} \S\otimes_\W\S&\mapright{\sigma\otimes\id}&\S\otimes\S\\
\lmapup{\gamma}&&\mapdown{\mu}\\
\S&\mapright{\kappa\circ\epsilon}&\S.\end{matrix}$$
Thus $\S$ is a weakly complete commutative symmetric Hopf algebra.

\msk  

Let us discuss its {\it primitive} and {\it grouplike}
elements:
 
\ssk\nin An element $\a=(a_n)_{n\in\N_0}\in\S$ is {\it primitive} 
if and only if $\Delta(\a)=\gamma(\a)$, that is, by (2) and (3), if and only if
$$ \sum_{m, n\in \N_0}(a_mX_1^m+ a_nX_2^n) = 
     \sum_{n\in\N_0}a_n(X_1+X_2)^n,\leqno(4)$$
if and only if  $n\ne 1\implies a_n=0$ if and only if $\a= t\.X$ for 
some $t\in \K$. Thus 
$$\P(\S)=\K\.X.\leqno({\rm PR})$$

\msk

An element $\a$ is {\it grouplike} if and only if it is nonzero and satisfies 
$\gamma(\a) =\a\otimes\a$ if and only if
$$ \sum_{n\in N_0}\, a_n(X_1+X_2)^n= \sum_{m,n\in\N_0}a_mX_1^m\.a_nX_2^n,\leqno(5)$$
which is the case if and only if 
$$(\forall n\in\N_0)\,a_n=\frac1{n!}.\leqno(6)$$
Thus
$$\G(\S)= \exp\K\.X=\begin{cases}\exp{\R\.X},\quad \mbox{if $\K=\R$},\\
(\exp \R\.X)(\exp \R i\.X),\quad \mbox{if $\K=\C$}.\end{cases}
\leqno({\rm GR})$$

This confirms the general result in \cite{dhtwo}, Theorem 6.15.

\bsk
\nin
It is now urgent that we precisely describe the {\it exponential function}
$\exp\colon\S\to \S^\times$:

For any field $\K$ we write $\K^\times$ for $\K\setminus\{0\}$.
We set $\I\defi\{\a: a_0=0\}$. Then $\I$ is the maximal ideal of $\S$
with $\S/\I\cong\K$. Notice that $\S=\K\.\1\oplus\I$ and $\exp(t\.\1+x)
=(\exp t)\.(\exp x)$ for $t\in\K$ and $x\in\I$. In particular,
$$\exp\S=(\exp \K)\.\exp(\1+\I),\mbox{  where} 
\exp \K=\begin{cases}\R_<=\{r\in\R:0<r\}& \mbox{if $\K=\R$},\\
                   \C^\times& \mbox{if $\K=\C$}.\end{cases}$$

\msk
\begin{Lemma} \label{logo} {\rm(i)} For $\a=(a_n)_{n\in \N}\in\S$ 
we have $\a\in\S^\times$ if and only if $a_0\ne 0$ and so  $\S^\times=\S\setminus\I$.  

{\rm(ii)} The function $\exp\colon\K\to\K^\times$ maps 
$\R$ bijectively onto $\R_<$ if $\K=\R$ and  $\C$ surjectively 
onto $\C^\times$ with kernel $2\pi i\Z$ if $\K=\C$.

The function 
$\exp\colon \I\to\1+\I$ is bijective with inverse $\log\colon \1+\I\to \I$.

{\rm(iii)} $\exp\colon 
\S=\K\.\1\oplus\I\to\S^\times=\K^\times\times(1+\I)$ 
is injective for $\K=\R$ and surjective for $\K=\C$.

{\rm(iv)} The function $\exp\colon \P(\S)\to\G(\S)$ is  bijective,  
the inverse function being the logarithm $\log$. 
\end{Lemma} 

\begin{Proof} (i) If $t\defi a_0\ne0$ then 
$\a= t\.(\1-Y)$ where $Y=(0,-a_1/t,a_2/t,\dots)$ and 
$(\1-Y)^{-1}=\1+Y+Y^2+\cdots $. So $\a$ is invertible.
On the other hand, if $a_0=0$, then $\a\b=(0,b_0,\dots)=(0,\dots)$ and so
$\a$ fails to be invertible. 

\msk

(ii) and (iii) were shown above.

\msk

(iv) This is immediate from the preceding, 
since $\P(\S)=\K\.X$ and  $\G(\S)=\exp(\P(\S))$.
\end{Proof}

We summarize our results on the power series algebra in one variable:

\begin{Proposition} \label{powers} {\rm (The power series algebra $\K[[X]]$)} 
{\rm(i)} The weakly complete power series 
algebra $\S\defi\K[[X]]=(\K^{\N_0},+,*)$ is a singly generated
weakly complete symmetric Hopf algebra  generated by the element $X$. 

\msk

\ssk\nin{\rm(ii)} $\S$ is a local weakly complete algebra with maximal ideal 
$$\I=\Big\{\sum_{n=1}^\infty a_nX^n: \a{=}(0,a_1,a_2,\dots){\in}\K^{\N_0}\Big\}
\mbox{ and }\S^\times{=}\S\setminus\I.$$
Further,
$$\S=\K\.\1\oplus \I\mbox{ and }(\forall t\in\K, x\in\I)\, \exp(t\.\1 + x)
=e^t\.\exp x.$$
\msk

\ssk\nin{\rm(iii)} $\exp\colon \I\to \1+\I$  has the inverse $\log$ 
and therefore implements
an isomorphism of pro-Lie groups $(\I,+)\cong (\1+\I,\.)$

\msk

\ssk\nin{\rm(iv)} The additive group $\P(\S)$ of primitive elements is $\K\cdot X$,
the multiplicative group $\G(A)$ of grouplike elements 
is $(\exp \K\.X,*)$. 
\end{Proposition}

Recall that on $\K$ (with $\T=\R/\Z$) we have
$$\exp\K=\begin{cases}\R_<=\{r\in\R:0<r\}\cong (\R,+),\mbox{ if $\K=\R$, and}\\
\C^\times\cong (\R\oplus\T,+),\mbox{ if $\K=\C$}.\end{cases}$$
By Proposition \ref{powers} we have the following isomorphisms of
abelian pro-Lie groups
$$(\S,+)=\K\.1\oplus \I\cong(\K,+)\times (\I,+)\cong (\K,+)^{\N_0},$$
and
$$(\S^\times,\.){=}(\S\setminus\I,\.){=}(\K^\times,\.)\.(\1+\I,\.)
{\cong}(\K^\times,\.)\times (\K,+)^\N,$$
$$ (e^\K,\.){\cong}
\begin{cases}(\R_<,\.),\mbox{ if $\K{=}\R$},\\
             (\C^\times,\cdot),\mbox{ if $\K{=}\C$}.\end{cases}$$
Note also that on the level of primitive and grouplike elements
we have simply
$$\P(\S)=\K\.X\cong (\K,+),\quad \G(\S)=\exp(\K\.X)\cong 
\begin{cases}(\R,+),\quad\mbox{$\K=\R$},\\
              \R\times\T,\quad\mbox{$\K=\C$}.\end{cases}.$$
Here one should keep  in mind
the example of the power series algebra over $\K$:
$$\S=\K\[X\]\quad(\cong \K^{\N_0}\quad \mbox{in}\ \W).$$ 

\bsk 

\subsection{The universal monothetic algebra} 

\ssk\nin We know from \cite{dhtwo} that there is a singly generated universal
weakly complete algebra $\K\<X\>=\UU_\K(\g)$, where $\g=\K$ is the 
one-dimensional Lie algebra.  At his point we shall also discuss 
the weakly complete symmetric Hopf algebra structure of $\K\<X\>$ 
For the following we refer to \cite{dhtwo}, Corollary 3.3ff.
The defining fact of $\K\<X\>$ is the following universal property:

\msk 
\nin
$\bullet$ \quad {\it For each weakly complete unital algebra $A$ in $\WA$ and each element
$a\in A$ there is a unique $\WA$-morphism $\phi\colon\K\<X\>\to A$ such that 
$\phi(X)=a$.}

\ssk
\nin We shall see that the internal structure of $\K\<X\>$ is more
complicated overall than one might expect initially.

\msk

It is clear that without loss of generality we may assume 
that $A$ is abelian.
By Theorem \ref{1.1}, we have $A\cong\lim_{I\in\JJ(A)} A/I$ 
and so $\bullet$ holds if
and only if it holds for all finite-dimensional commutative algebras $A$.
However, this universal property is satisfied exactly by 
the weakly complete
algebra $\lim_{J\in\I(\K[x])}\K[x]/J$ for the polynomial 
ring $\K[x]$ in one variable
$x$ over $\K$ and the filter basis of all of its ideals 
$\I(\K[x])$. Since $\K[x]$
is a principal ideal domain, every $J\in \I(\K[x])$ is of the form 
$J=(f)=f\K[x]$ for some polynomial $f\in\K[x]$. 
We may assume
$$\K\<X\>=\lim_{J\in\I(\K[x])}\frac{\K[x]}{J}
=\lim_{f\in\K[x]}\frac{\K[x]}{(f)}
\subseteq \prod_{f\in\K[x]}\frac{\K[x]}{f\K[x]},\leqno(1)$$
generated by $X\defi(x+f\K[x])_{f\in\K[x]}\in\K\<X\>$.
(See also \cite{dhtwo}, Lemma 3.4.)

\ssk

We let $\PP=\PP_\K$ denote the set  of the 
irreducible polynomials $p$ over $\K$ with leading coefficient $1$ from
the polynomial ring $\K[x]$. Then $f\in\K[x]$, 
by the Chinese Remainder Theorem, is of the form
$$ f=t\.\prod_{p\in\PP}p^{k_p}\mbox{ some $t\in\K$ and $(k_p)_{p\in\PP}\in{\N_0}^{(\N_0)}$},$$
where ${\N_0}^{(\N_0)}$ denotes the set of all families of nonnegative integers vanishing
with the exception of indices $p$ from a finite subset of  $\PP$.
Then for all $f\in \K[x]$ we have.
$$\frac{\K[x]}{(f)}\cong \prod_{p\in\PP}\frac{\K[x]}{(p^{k_p})}.
\leqno(2)$$
Now $\{\K[x]/p^n\K[x]:n\in\N\}$ is a projective system and we introduce the
notation  
$$\S_p\defi\lim_{n\in\N} \frac{\K[x]}{p^n\K[x]}
\leqno(3)$$
which  is 
is a weakly complete commutative algebra generated by
$$X_p\defi(x +p^n\K[x])_{n\in\N}\in \lim_{n\in\N}\frac{\K[x]}{p^n\K[x]}.$$
In $\S_p$ we have a maximal ideal $\I_p:=X_p\S_p$ so that 
$\S_p/\I_p\cong\K$ and 
$$ \S_p=\K\.\1\oplus \I_p.$$
We conclude
$$\K\<X\>=\prod_{p\in\PP_\K}\S_p,\quad X=(X_p)_{p\in\PP_\K}.
\leqno(4)$$
In the interest of brevity again, we shall also write $\SS$
in the place of $\K\<X\>$. For easy reference we summarize 
the preceding  discussion
in the following lemma:
 
\begin{Lemma} \label{ffactor}{\rm(i)} We have
 $\SS=\prod_{p\in\PP_\K}\S_p$.
For each $p\in\PP$ the algebra $\S_p$ is
generated algebraically and topologically by $X_p$,
and $\SS$ is generated algebraically and topologically by 
$X=(X_p)_{p\in\PP}\in\SS$.

{\rm(ii)} For each weakly complete algebra $A$ and each $a\in A$ there is a
unique $\WA$-morphism     $\phi_a\colon\SS\to A$ such that 
$\phi_a(X)=a$.
\end{Lemma}

\bsk

The remainder of this section is devoted to a clarification
of  the structure of $\S_p$ defined in (3). 

\begin{Lemma} \label{insert} If $p\in\PP$ is of degree 1, then
$$\S_p\cong \K\[X_p\].\leqno(5)$$
\end{Lemma}
\begin{Proof} If $c\in\K$ and $p=x-c$, abbreviate $R:=\K[x]\cong \K[x-c]$.
Then $x\mapsto x-c$ induces an automorphism of $R$ and
$\S_p=\lim_{n\in\N}\frac{R}{(x-c)^nR}\cong
       \lim_{n\in\N}\frac{R}{x^nR}\cong \K\[X_p\]$.
\end{Proof}

Since for $\K=\C$ every $p\in\PP_\C$ is of degree 1, we know
that in this case $\S_p\cong \C[[X_p]]$ for all $p\in\P$.

\msk

Now we assume $\K=\R\subseteq\C$. Then there are two cases:

\nin (a) $p=x-r$ for  $r\in \R$. 
Then  $$\S_p\cong\R[[X_p]].\leqno(5a)$$

\msk

\nin (b) There is a $c\in \C\setminus \R$   such that 
$$p(x)=(x-c)(x-\overline c)=x^2-(c+\overline c)x +c\overline c=
x^2-2{\rm Re}(c)x+|c|^2,\quad p\in\PP_\R.$$
In this case we write 
$p_1=x-c$, Im\ $c>0$, and $p_2=x-\overline c$,\ $p_n\in\PP_\C$, $n=1,2$.

\begin{Lemma} \label{linus1} In the case of {\rm (b)} above, the real algebra 
$\frac{\R[x]}{p^n\R[x]}$ is isomorphic to the
real algebra underlying $\frac{\C[x]}{p_1^n\C[x]}$.
\end{Lemma}
\begin{Proof} We use the abbreviations $R_n=\frac{\R[x]}{p^n\R[x]}$,
$C_n=\frac{\C[x]}{p^n\C[x]}$, and
$C_{kn}=\frac{\C[x]}{p_k^n\C[x]}$, $k=1,2$.
By the Chinese Remainder Theorem we have an isomorphism
$$\rho\colon C_n\to C_{1n}\times C_{2n}\leqno({\rm i})$$  
such that $\rho(u+p^n\C[x])=(u+p_1^n\C[x],u+p_2^n\C[x])$. 
The real algebra  underlying the right hand side of (i)
has an involution $\sigma$ defined by
$$\sigma(u+p_1^n\C[x],v+p_2^n\C[x])
=(\overline v+p_1^n\C[x],\overline u+p_2^n\C[x]),\leqno({\rm ii})$$
such that the elements of the
 real fixed point algebra $F$ of $\sigma$ are the
elements $(u+p_1^n\C[x],\overline u+p_2^n\C[x])$ with
$u\in\C[x]$. The restriction of the projection
%$$\frac{\C[x]}{p_1^n\C[x]}\times \frac{\C[x]}{p_2^n\C[x]}
$C_{1n}\times C_{2n}\to C_{1n}$ to $F$ is an isomorphism.
Thus $F\cong C_{1n}$ as real algebra. In particular,
 $\dim_\R F =\dim_\R C_{1n}= 2n =\dim_\R R_n$. Thus
the injection $R_n\to F$ via $\rho$ is in fact surjective.
Hence $\frac{\R[x]}{p^n\R[x]}=R_n\cong F=\frac{\C[x]}{p_1^n\C[x]}$
as real algebras.
\end{Proof}
 
As a consequence we conclude that 
$$\S_p=\lim_{n\in\N}\frac{\R[x]}{p^n\R[x]}
\cong\lim_{n\in\N}\frac{\C[x]}{p_1^n\C[x]}$$
and thus
$$\S_p\cong \C\[X_{p_1}\]\mbox{ as real algebras}.\leqno(5b)$$ 

\msk

For $\K=\R$ or $\K=\C$, there is an injection $p\mapsto c_p:\PP_\K\to \C$, 
where 
$$ p(x)=\begin{cases} x-c_p,\quad\mbox{if $\deg(p)=1$,}\\
 x-c_{p_1},\ {\rm Im\ }c_{p_1}>0,\quad \mbox{if $\deg(p)=2$.}
\end{cases}\leqno{(6')}$$
If $\K=\C$, then we are in the first case and $c_p$ ranges through all of
$\C$. If $\K=\R$, then both cases occur, and in the first case 
$c_p$ ranges through $\K=\R$ and in the second case $c_{p_1}$ ranges
through the open upper complex half-plane.% $\C^\uparrow$.

\msk

The different cases now sum up to the following statement:

\begin{Lemma} \label{twocases} 
$$ \K\<X\>{=}\prod_{p\in\PP_\K}\S_p{\cong}
\begin{cases} \prod_{p\in\PP_\R,\deg p=1}\R\[X_p\]
        \times\prod_{p\in\PP_\R,\deg p=2}\C\[X_{p_1}\],\,\mbox{if $\K=\R$},\\
              \prod_{p\in\PP_\C}\C\[X_p\],\quad \mbox{if $\K=\C$},\end{cases}
\leqno{(6'')}$$
where all algebras in the top line are real algebras. \end{Lemma}

\msk 

\nin For $p\in\PP_\K$ we write
$$\K_p=\begin{cases}
\R,\quad \mbox{if $\K=\R$ and $\deg p=1$,}\\
\C,\quad \mbox{if either $\K=\R$ and $\deg p=2$, or $\K=\C$.}
        \end{cases}\leqno({\bf*})$$
Moreover, elements in $\S_p$ with $p\in \PP_\K$ we denote by
$$\a_p=\sum_{n\in\N_0}a_{np}{X_p}^n, \quad a_{np}\in\K_p.$$
Finally  we have 
$$\SS\defi\K\<X\>=\prod_{p\in \PP_\K}\S_p
=\{\aa\defi(\a_p)_{p\in\PP_\K}: \a_p\in \S_p\}, \leqno(7)$$
a weakly complete commutative symmetric Hopf algebra, with
componentwise operations and co-operations. 
By Lemma \ref{ffactor}(i), the weakly complete
algebra $\SS$ is the algebraically and topologically singly generated  weakly
complete algebra  with generator 
$$X\defi(X_p)_{p\in \PP_\K}\in \prod_{p\in\PP_\K}\S_p,\leqno(7X) $$

\bsk 

\nin For each $p\in\PP_\K$ we define $\I_p\subseteq\S_p$ 
to be the maximal ideal of $\S_p$, where $\exp\colon \I_p\to 1+\I_p$ has 
the inverse $\log\colon1+\I_p\to\I_p$. 
In view of Lemma \ref{logo} we observe the following fact:

\begin{Remark} \label{expo} The exponential function $\exp\colon\SS\to\SS^\times$
is given componentwise for $\aa=(\a_p)_{p\in\PP}$ as
$$\exp\aa=(\exp \a_p)_{p\in\PP}.$$
The exponential function  $\exp\colon \S_p\to\S_p$ is surjective 
if either $\K=\R$ and $\deg p=2$, or $\K=\C$, 
and it is injective if $\K=\R$ and $\deg p=1$. 
\end{Remark} 

\bsk

We now aim to discuss the restriction of the exponential function to
the set of primitive elements. First
recall that  on $\SS=\K\<X\>$ we have the  operations
$$\aa+\bb= (\a_p+\b_p)_{p\in\PP}\mbox{ and }
  \aa*\bb= (\a_p*\b_p)_{p\in\PP}
=\bigg(\sum_{n\in\N_0}
\Big(\sum\nolimits_{k+m=n}a_{kp}b_{mp}\Big){X_p}^n\bigg)_{p\in\PP}.$$
In $\SS\otimes_\W\SS$, with $\PP=\PP_\K$  we write
$$X_{1p}:=(X_p)_{p\in \PP}\otimes \1=(X_p\otimes \1)_{p\in\PP}\in\SS\otimes_\W \SS,$$
and
$$X_{2p}:=\1\otimes (X_p)_{p\in \PP}=(\1\otimes X_p)_{p\in\PP}\in\SS\otimes_\W\SS,$$
Again we may consider $\SS\otimes_\W\SS$ as weakly complete power series algebra
with commuting variables $X_{1p}$ and $X_{2p}$ with $p$ ranging through $\PP$.

\bsk

The identity and coidentity $\epsilon\colon \K\to \SS$, 
$\kappa\colon \SS\to\K$ (augmentation), and symmetry
$\sigma\colon \SS\to \SS$ are straightforward from the
respective operations in $\S$,
but let us also consider  the diagonal vector space
morphism $\Delta$ and the algebra comultiplication $\gamma$:
$$\Delta,\gamma: \SS\to\SS\otimes_\W\SS\mbox{\quad defined as follows:}$$
For $\aa=\sum_{(n,p)\in\N_0\times\PP} a_{np}{X_p}^n$  we have
$$ \Delta(\aa){:=}\aa\otimes\1+\1\otimes\aa
  {=}\sum_{m,n\in\N_0,f,g\in\PP}(a_{mp}{X_{1p}}^m{+}a_{ng}{X_{2g}}^n)\leqno(4)$$
and
$$ \gamma(\aa)=\sum_{n\in\N_0,p\in\PP} a_{np}(X_{1p}+X_{2p})^n.\leqno(5)$$
We have the commutative diagram
$$\begin{matrix} 
\SS\otimes_\W\SS&\mapright{\sigma\otimes_\W\id}&\SS\otimes_\W\SS\\
\lmapup{\gamma}&&\mapdown{\mu}\\
\SS&\mapright{\kappa\circ\epsilon}&\SS\end{matrix}$$
identifying  $\SS$ as  weakly complete symmetric Hopf algebra allowing
us now to 
turn  to the determination of the primitive  and grouplike elements 
of $\SS$. 
Indeed an element $\aa$ is {\it primitive} in $\SS$ 
if $\gamma(\aa)=\Delta(\aa)$, that is, if and only if
$$\sum_{m,n\in\N_0,\ p\in \PP}(a_{mp}{X_{1p}}^m+a_{nf}{X_{2p}}^n)
=\sum_{n\in\N_0,\ p\in\PP}a_{np}(X_{1p}+X_{2p})^n$$
if and only if
 $$(\forall p\in\PP)\, n\ne1\implies a_{np}=0$$
if and only if
$$(\forall p\in\PP)(\exists t_p\in \K_p)\, \aa=(t_p X_p)_{p\in\PP}.$$
Thus we have
$$\P(\SS)= \prod_{p\in\PP_\K}\K_p\cdot X_p.\leqno({\rm PR})$$
\bsk
On the other hand, 
an element $\aa$ is {\it grouplike} if it is nonzero and satisfies
$\gamma(\aa)=\aa\otimes\aa$, that is,
$$\sum_{n\in\N_0,p\in\PP} a_{np}(X_{1p}+X_{2p})^n=
\sum_{m,n\in\N_0, p\in\PP}a_{m,p}X_1^m\.a_{ng}X_2^n,$$
which is the case if and only if
$$(\forall (n,p)\in\N_0\times\PP)\  a_{np}=\frac1{n!}.$$
Thus we have
$$\G(\SS)=\prod_{p\in\PP_\K}\exp(\K_p\.X_p).\leqno({\rm GR})$$
 
\begin{Definition} \label{unimono} The 
algebraically and topologically singly generated
universal weakly complete algebra $\SS=\K\<X\>$ 
is called the {\it universal monothetic algebra}. 
\end{Definition}

\nin Lemma \ref{ffactor} (i) justifies
 the name.
Recall from Proposition \ref{powers} that for each $p\in\PP_\K$ the algebra
$\S_p$ is a local weakly complete algebra
 with a  maximal ideal $\I_p$ and that 
${\S_p}^\times=\K_p^{\times}\.(\S_p\setminus \I_p)$ where $\K_p$ was defined
in $({\bf*})$. 
Now we are prepared to summarize the structure theorem for $\SS$:
\begin{Theorem} \label{monoth} 
{\rm (The universal monothetic algebra $\K\<X\>$)} 

\nin{\rm(i)} The universal monothetic 
algebra 
$$\SS\defi\K\<X\>=\prod_{p\in\PP_\K}\S_p$$ is a 
weakly complete symmetric Hopf algebra  generated by the element 
$X=(X_p)_{p\in\PP_\K}$. 

\msk

\ssk\nin{\rm(ii)} The group of units $\SS^\times$ is dense in $\SS$, where  
$\SS^\times=\prod_{p\in\PP_\K}{\S_p}^\times=$
$$\Big\{\big(\sum\nolimits_{n\in\N_0} a_{np}{X_p}^n\big)_{p\in\PP_\K}:  
\ a_{np}\in \K_p \mbox{ and }
 a_{0,p}{\ne}0\Big\}$$
$$=\prod_{p\in\PP_\K}\K^\times_p(\S_p\setminus \I_p).\leqno(\#)$$

\msk

\ssk\nin{\rm(iii)} The exponential function $\exp\colon \SS\to\SS^\times$ operates
componentwise on $\prod_{p\in\PP_\K}\S_p$ and induces an isomorphism of 
of topological groups 
 $\prod_{p\in\PP_K}\I_p\to \prod_{p\in\PP_\K}(1+\I_p)$, 
whose inverse is given by the componentwise logarithm.

\msk

\ssk\nin{\rm(iv)} The additive group $\P(\SS)$ of primitive elements is 
$$\P(\SS)=\prod_{p\in\PP_\K} \K_p\cdot X_p\subseteq
 \prod_{p\in\PP_\K}\S_p.$$ 
In particular, the element $X$ is primitive.
The image of $\g$ in $\K\<X\>$ is $\K\.X$, i.e.\
 $\g\subset \P(\UU(\g))$.
\msk

\nin  Let $\T=\R/\Z$ denote the additive circle group again and $\c$
the cardinality $2^{\aleph_0}$ of the continuum.
The multiplicative group $\G(\SS)$ of grouplike elements is 
$$\G(\SS)=\prod_{p\in\PP_\K} \exp(\K_p\.X_p)\cong 
\begin{cases}\R^{\PP_\R}\oplus
\T^{\{p\in\PP_\R:\ \deg p=2\}},\,\mbox{if $\K=\R$},\\
 (\R\oplus\T)^{\PP_\C},\, \mbox{if $\K=\C$}\end{cases}
\hskip-15pt\mbox{$\Bigg\}$}
\cong(\C^\times)^\c.$$
The exponential function $\exp\colon\P(\SS)\to \G(\SS)$ is  a quotient
morphism of topological abelian groups onto its image.
\end{Theorem}

\nin In particular we derive the following (with $\c=2^{\aleph_0}$):

\begin{Corollary}\label{grouplike-conn} The abelian pro-Lie group
$\G(\K\<X\>)$ is connected and, for $\K=\R$, 
is isomorphic to $\R^\c\oplus \T^\c\cong (\R\times\T)^\c\cong(\C^\times)^\c$. 
\end{Corollary}

We do not know whether in general the pro-Lie group $\G(\UU(\g))$
is connected.
 Theorem \ref{monoth} (iv)) shows that $\g\subseteq A\defi\UU_\K(\g)$ is
considerably smaller than $\P(A)$. The discrepancy between 
$\g$ and  $\P(A)$ arises
in the detailed  description of 
the universal monothetic algebra 
$\K\<X\>$. The  origin  of this  complication 
is the Galois theory of the polynomial ring $\K[x]$.

\msk

\subsection{Comments }

\ssk\nin We discussed extensively the ``smallest possible'' 
nontrivial weakly
complete enveloping algebra $\UU_\K(\g)$, namely, the one arising
for $\dim\g=1$.
 Any {\it abelian} profinite-dimensional
Lie algebra $\g$ is isomorphic to $\K^J$ for some set $J$.
We have a pair of adjoint functors between 
the category $\W$ of weakly complete
 vector spaces over $\K$
and the category $\WAC$ of weakly complete 
commutative unital algebras, namely,
the functor $A\mapsto |A|: \WAC\to \W$ 
assigning to a weakly complete commutative algebra its
underlying weakly complete vector space and 
 $\g\mapsto \UU_\K(\g):\W\to \WAC$, the restriction 
 of the universal enveloping functor. 
Then $\UU_\K|\W$ is left adjoint to $|\cdot|$,
and therefore it  preserves colimits.   For
a finite set $J$ of $n$ elements we note that in the category $\W$
we have
$\g=\g_1\oplus\cdots\oplus\g_n$ with $\g_k\cong\K$ for $k=1,\dots,n$
and so $\g$ is the coproduct of $n$ cofactors of dimension 1. Accordingly, 
in the category $\WAC$ we observe
$$\UU_\K(\g)\cong \coprod_{j=1}^n\UU_\K(\g_j),\quad \UU_K(\g_j)
\cong\SS\cong\S^{\PP_\K}.\leqno(1)$$ 
To the extent that finite coproducts in the category $\WAC$
are understood, one knows $\UU_\K(\g)$ for finite
dimensional abelian Lie algebras $\g$.
\msk

Let now $\g=\K^J$ in $\W$ for some nonempty  set $J$. 
If $J$ is finite, then $\g$ is a finite coproduct and the dual of a finite product.
If $J$ is infinite, then let $\fin J$ denote  the directed
family of finite subsets $F\subseteq J$ and recall
$$\g=\K^J\cong\lim_{F\in\fin J}\K^F.$$
We may then apply Theorem \ref{lim2} and deduce
$$\UU_\K(\g)\cong \lim_{F\in\fin J}\UU_\K(\K^F),\leqno(2)$$
where $\UU_\K(\K^F)$ is known by (1) if we know finite
coproducts in the category $\WAC$.

\msk
 
In Theorem \ref{completion} in the appendix we argued that for any 
profinite-dimensional Lie algebra $\g$ with the underlying weakly complete
vector space $|\g|$ there is a canononical quotient morphism of weakly complete 
unital algebras $q_\g\colon \TT(|\g|)\to\UU(\g)$ from the weakly
complete tensor algebra $\TT(|\g|)$ onto $\UU(\g)$. In the present
situation  of abelian Lie algebras we may write $|\g|=\g$ and
conclude that for each set $J$ and for $\g=\K^J$ we have
 a natural quotient morphism of
weakly complete algebras 
$$ q_\g\colon \TT_\K(\g)\to \UU_\K(\g) \leqno(3)$$
whose kernel is the closed ideal generated in $\TT_\K(\g)$ by the
elements $xy-yx$, $x,y\in\g\subseteq\TT(\g)$. 
The structure theory of $\SS$ for
the case $\g=\K$ shows that the presentation (3) conceals more than
it reveals. Indeed, the universal property of $\UU_\K(\K^J)$ yields
for $\g=\K^J$ a surjective morphism of weakly complete unital algebras
$$A\defi\UU_\K(\g)=\UU_\K(\K^J)\to \UU_\K(\K)^J\cong\SS^J
\quad(\cong \K^{\N_0\times \PP_\K\times J}\mbox{ in }\W).\leqno (4)$$

\msk
Since $\P(\SS^J)=\P(\SS)^J\cong \K^{\PP\times J}$
the quotient morphism in (4) shows that the vector space of
primitive elements of $\SS^J$ is considerably larger than $\g=\K^J$.
This information indicates that for the weakly complete
symmetric Hopf algebra $A:=\UU(\g)$, the subspace $\P(A)$ of primitive elements is
likely to be  large by comparison with $\g$.
Clearly $\G(\SS^J)$ is a {\it connected} abelian group whose
structure is known by Corollary  \ref{grouplike-conn}. 
The simplest example along this line is the power series 
Hopf algebra $\K\[X\]$ (cf.\ Proposition \ref{powers}). Accordingly,
one expects the group $\G(A)$ to be considerable.
\bsk

However, some caution is in order:

\begin{Example} \label{dense-exp}  {\rm (i)} Let $A=\R[\hat\Q]$ be the real
group algebra of 
$\hat\Q$,  the universal solenoidal compact abelian group. Then
$\G(A)=\hat\Q\subseteq A$, while $\P(A)=\L(\hat\Q)\cong\R$. Then
$\exp_{\hat\Q}\colon\P(A)\to\G(A)$ is a morphism of locally compact
abelian  groups with
 a dense image, but it is {\it not} surjective.

\msk
{\rm(ii)} If we take  $A=\R[\Z_p]$ for a prime number $p$, where $\Z_p$
is the additive group of the $p$-adic integers,  then
$\G(A)\cong \Z_p$ and  $\P(A)=\L(\G(A))=\{0\}$, since $\Z_p$ is totally
disconnected. 
The exponential map  $\exp_A\colon \P(A)\to \G(A)$ is the zero morphism.
\end{Example}
 
These examples show that even on the abelian level, the weakly complete
symmetric Hopf algebra structure of the weakly complete enveloping algebras
and that of the weakly complete group algebras behave rather differently.
Yet they are related in a natural way as we shall observe in the following
section.

\section{Enveloping Algebras Versus Group Algebras}

The class of compact groups and their Lie algebras are
 distinguished domains for which the 
relationship between weakly complete enveloping algebras and weakly complete
group algebras is particularly lucid. Hence we focus on these classes.

\subsection{The case of compact groups}
A particularly appropriate situation is that of a {\it compact} topological
group $G$. Our level of information regarding the associated {\it real} group
algebra is particularly advanced in that situation. Indeed recall that for a
compact group we may naturally identify $G$ with the group of 
grouplike elements of $\R[G]$ (cf.\ \cite{dhtwo}, Theorems 8.7,
8.9 and 8.12), and we may further identify  $\g\defi\L(G)$  with 
the pro-Lie algebra  $\P(\R[G])$ of
primitive elements. (Cf.\  also Theorem \ref{expone} in the Appendix.)
We  may also assume that the Lie algebra $\g$ of $G$ is contained in 
the set $\P(\UU_\R(\g))$ of primitive
elements of $\UU_\R(\g)$.

\begin{Theorem} \label{8.7} {\rm(i)}  Let $G$ be a 
compact group and $\g$ its Lie algebra. Then there
is a natural morphism of weakly complete algebras 
$\omega_G\colon\UU_\R(\g)\to \R[G]$ fixing the elements
of $\g$ elementwise.

{\rm(ii)}  The image of $\omega_G$ is the closed subalgebra 
$\R[G_0]$ of $\R[G]$.

{\rm(iii)}  The pro-Lie group $\G(\UU_\R(\g))$ is mapped into 
$G_0=\G(\R[G_0])\subseteq\R[G]$. The connected pro-Lie group 
$\G(\UU_\R(\g))_0$ maps epimorphically to  $G_0$ and $\P(\UU_\R(\g))$
maps surjectively 
onto $\P(\R[G])=\g$.
\end{Theorem}

\begin{Proof} (i) follows at once from the universal property of $\UU$. 

(ii) As a morphism
of weakly complete Hopf algebras, $\omega_G$ has a closed image which is
generated as a weakly complete subalgebra by $\g$ which is
$\R[G_0]$ by Corollary 3.3 (ii) of \cite{hofkra}. 

(iii) The morphism $\omega_G$ of weakly complete Hopf algebras maps
grouplike elements to grouplike elements, 
whence we have the commutative diagram
$$\begin{matrix}
\g\subseteq\P(\UU_\R(\g))&\mapright{\P(\omega_G)}&
       \P(\R(G)){=}\g\\
\hskip40pt\lmapdown{\exp_{\G(\UU_\R(\g))}} &&\mapdown{\exp_G}\\
\hfill\G(\UU_\R(\g))&\lmapright{\G(\omega_G)}&
                         \G(\R[G])=G.
\end{matrix}$$
Since $\P(\omega_G)$ is a retraction and the image of $\exp _G$
topologically generates $G_0$,  the image
of $\G(\omega_G)\circ \exp_{\UU_\R(\g)}$ topologically generates
$G_0$. Since the image of the exponential function of the pro-Lie group 
$\G(\UU_\R(\L(G)))$ generates topologically its identity
component,  $\G(\omega_G)$ maps this identity component onto $G_0$.

Since $\g\subseteq \P(\UU_\R(\g))$, and since also any morphism
of Hopf algebras maps a primitive element onto a primitive
element we know $\omega_G(\P(\UU_\R(\g)))=\P(\R[G])$.\end{Proof}

\msk

\medskip

\noindent The following overview of the situation may be  helpful:

\vskip-25pt
$$
\begin{matrix}%
   &&&& \R[G]\\ 
   &&&& \Big|\\
\UU_\R(\g)&\mapright{\omega_G,{\rm onto}}&\R[G_0]&&|\\
 \Big| &&\Big|&&\Big|\\
\G(\UU_\R(\g))&\mapright{}&\big|&\mapright{}&\G(\R[G])=G\\
\Big|&&\Big|&&\Big|\\
\G(\UU_\R(\g))_0&\mapright{\rm onto}&\G(\R[G_0]){=}G_0&
                              = &G_0\\
\lmapup{\exp_{\G(\UU_\R(\g))}}&&\mapup{\exp_{\G(\R[G])}}&=&\mapup{\exp_G}\\
\g{\subseteq}\P(\UU_\R(\g))&\lmapright{\rm retract}&\P(\R[G_0]){=}\P(\R[G])&=&\g.\\
\end{matrix}\leqno{(D)}$$

\msk 

\begin{Example} \label{ex1} Let $\g$ be a compact semisimple Lie algebra. Then $\g=\L(G)$
for the compact projective group $G=\Pr(\g)$. In this case, $G=\Gamma(\g)$, and we have
a commutative diagram

$$\begin{matrix}
\UU(\g)&\mapright{\omega_G,surjective}&\R[G]\\
\Big|&&\Big|\\
\G(\UU(\g))_0&\mapright{\rm retract}&G=\G(\R[G])\\
\lmapup{\exp_\G}&&\mapup{\exp_G}\\
 \g\subseteq \P(\UU(\g))&\lmapright{\rm retract}&\g.\\
\end{matrix}\leqno(D_1)$$

\end{Example}
We do not precisely know what $\G(\UU(\g))$ and $\P(\UU(\g))$
are even if $\g={\rm so}(3)$ in which case $\Gamma^*(\g)\cong {\rm SU}(2)$.
Still, in this case $\exp_{\Gamma(\g)}\colon \g\to\Gamma(\g)$
is surjective (cf.\ \cite{compbook}, Theorems 6.30,  
9.19(ii)  and Theorem 9.32(ii)).  

 The group $\G(\UU(\g))$ of grouplike elements
of $\UU(\g)$ is a semidirect product of some unknown closed normal
subgroup  $N$ by $G$. From the content of Diagram $(D_1)$ we do not
know anything about $N$. 

\msk

The following example is the opposite to the preceding one:
  
\begin{Example} \label{ex2} Let $\g=\R^X$ for some set $X$.
Then $G=\Pr(\g)=(\hat\Q)^X$ and $\Gamma(\g)=\Gamma^*(G)=\R^X$.
\end{Example}

In our discussion of abelian profinite-dimensional Lie algebras $\g$
we have obtained more information on $\UU(\g)$. Here we have our standard 
diagram:

$$\begin{matrix}
      \UU(\g)& \mapright{\omega_\g,\rm surjective} &\R[G]{\subset}\C^{\Q^{(X)}}\\
      \Big|  &             &\Big|\hfill\\
\G(\UU(\g))_0& \mapright{\rm onto} &G=(\hat\Q)^X\hfill\\
\lmapup{\exp_\G}& &\mapup{\exp_G}\hfill\\
\g\subseteq\P(\UU(\g))&\lmapright{\rm retract}&\g.\hfill\\
\end{matrix}\leqno(D_2)$$

\bsk

Recall that a finite-dimensional real Lie algebra $\g$ is called ``compact'' 
if it is isomorphic
to the Lie algebra of a compact group (apologetically defined in \cite{compbook}
Definition 6.1 in that fashion). 
We now expand this definition to read as follows:

\begin{Definition} \label{apology} A  Lie algebra is called
{\it compact} if it is profinite-dimensional and is 
isomorphic to the Lie algebra of a compact group.
\end{Definition}

We know a real Lie algebra to be {\it compact}
if and only if there
exists a set $X$ and a family $\SI$ of compact simple 
Lie algebras  $\s$ such that $\g\cong \R^X\times \prod\SI$,
where we wrote $\prod\SI$ for $\prod_{\s\in\SI}\s$. Now  from \cite{compbook},
Theorem 9.76 we obtain the following statement:

\begin{Theorem} \label{lie3c} {\rm(Sophus Lie's Third Principal Theorem
for Compact Lie Algebras)}  
For every compact
real Lie algebra $\g$ there is a projective connected compact group $\Pr(\g)$ 
whose Lie algebra $\L(\Pr(\g))$ is (isomorphic to) $\g$.
\end{Theorem} 

Every compact connected group $G$ with $\L(G)\cong\g$ is 
a quotient of $\Pr(\g)$.
modulo some central $0$-dimensional subgroup.
For details see \cite{compbook}, discussion following Lemma 9.72, notably
Theorem 9.76 and Theorem 9.76bis. For the abelian case see \cite{compbook},
Theorem 8.78ff.
Notice that for a compact  Lie algebra $\g$ the projective compact connected
group $\Pr(\g)$ is simply connected if and only if $\g$ is semisimple. 
By contrast, if
$\g=\R^X$ for some set then $\Pr(\g)=(\hat\Q)^X$ 
(see \cite{compbook}, Proposition 8.81),
a compact connected abelian group 
that fails to be simply connected while $\pi_1(\Pr(\g))=\{0\}$ 
(see \cite{compbook}, Theorem 8.62).

It is very important here to distinguish between the prosimply 
connected pro-Lie
group $\Gamma(\g)$ and, in the case of a compact  
Lie algebra $\g$, the projective compact group $\Pr(\g)$.

\bsk

The present concept of weakly complete enveloping algebras 
now belongs to the 
circle of ideas of Lie's Third Fundamental Theorem.

\msk

 Let $\g$ be a profinite-dimensional Lie algebra over $\R$.
By Theorem \ref{8.5} above, $\Gamma^*(\g)\subseteq\G(\UU(\g))$ 
is a pro-Lie group whose Lie algebra is $\g$ and the exponential 
function $\exp\colon \li{\UU(\g)}\to \UU(\g)^\times$ of $\UU(\g)$ 
induces the exponential function 
$$\exp: \L(\Gamma^*(\g))=\g \to \Gamma^*(\g).$$

If $G=G_0=\Gamma(\g)$ embeds into its weakly complete group algebra $\R[G]$,
then the diagram $(D)$ above shows that 
$\alpha_\g$ is an isomorphism. 

We summarize for $\K=\R$, 
recalling that we consider $\g$ as a Lie subalgebra of 
$\P(\UU(\g))\subseteq \UU(\g)$. Indeed,
in the context of Lie's Third Fundamental Theorem, 
there are two basic pro-Lie groups $\Gamma(\g)$
and $\Gamma^*(\g)$ attached, and, in the case of a 
compact profinite-dimensional Lie algebra $\g$,
a third one, $\Pr(\g)$, and for these we have:

\begin{Theorem} \label{U-projective} Let $\g$ be a 
profinite-dimensional real Lie algebra. Then
the pro-Lie group of  grouplike elements in the weakly 
complete enveloping algebra $\UU(\g)$
contains the pro-Lie group $\Gamma^*(\g)$, having $\g$ as 
Lie algebra with the 
exponential function $\exp\colon \g\to \Gamma^*(\g)$ 
induced by the exponential function
of $\UU(\g)$. There is a natural quotient morphism 
$\alpha_\g\colon \Gamma(\g)\to\Gamma^*(\g)$.
If the natural morphism
$\Gamma(\g)\to\R[\Gamma(\g)]$ of $\Gamma(\g)$ into 
its weakly complete group algebra is 
an embedding as is the case if $\g$ is a compact  Lie algebra, then
$\alpha_\g$ is an isomorphism, and in the latter case, 
there is a natural injective morphism 
$\Gamma(\g)\to \Pr(\g)$ with dense image.
\end{Theorem}

\bsk

\section{Appendix: The Category Theoretical Background}
 
For a category $\TA$ of topological algebraic structures---in 
the simplest case the category $\W$ of weakly complete vector spaces,
and for the category $\AA$ of weakly complete associative unital algebras,
we shall repeatedly discuss an adjoint pair of functors $R\colon \AA\to \TA$ 
and $L\colon \TA\to \AA$. As an example on the simplest level,
in the case of $\TA=\W$, for a weakly complete algebra 
$A$, the $\W$-object $R(A)$ will simply be the weakly complete vector space
underlying $A$, while for a  weakly complete vector space $W$, the
weakly complete algebra $L(W)$ will be the weakly complete
tensor algebra of $W$ in the category $\W$.

\bsk

\subsection{Limits and topologically dense subcategories}

\msk

Since the category $\AA$ of weakly complete associative  unital algebras
is at the focus of our considerations, let us point to one important
property of the objects in this category, which was expressed
in Appendix 7 of \cite{compbook}, Theorem A7.34.

\msk

\begin{Theorem} \label{1.1} For every weakly complete  
unital topological $\K$-algebra $A$, the set $\JJ(A)$ of closed
two-sided ideals $I$ with finite-dimensional quotient
algebras $A/I$ is a filterbasis converging to $0$ in $A$, and $A$ is
(naturally isomorphic to) the  projective limit 
$\lim_{I\in\JJ(A)} A/I$ of these 
{\em finite-dimensional} unital quotient algebras.
\end{Theorem}

\bsk   

This theorem says that any weakly complete unital associative
algebra ``is approximated by finite-dimensional $\K$-algebras''.
Let us briefly recall our approach to projective limits in a
category $\cA$. Each directed set $J$ is a category with 
the elements of $J$ as
objects and for each pair $(j,k)$ satisfying $j\le k$ an 
arrow ($J$-morphism) $k\to j$. A {\it projective (or inverse) system} is a
functor $J\to \cA$, usually written $j\mapsto A_j$ and $(k\to j)
\mapsto (f_{jk}\colon A_k\to A_j)$. The {\it projective limit of
this system} is an object $\lim_{j\in J}A_j$ together with a family
of morphisms $f_k\colon \lim_{j\in J}A_j\to A_k$, $k\in  J$
such that $f_k=f_{kn}f_n$ for all arrows $n\to k$.
The limit has the universal property that for any system
of morphisms $\phi_k\colon A\to A_k$ of $\cA$-morphisms satisfying
$\phi_k=f_{kn}\phi_n$ for all arrows $n\to k$ there is a 
{\it unique} morphism $\phi\colon A\to\lim_{j\in J}A_j$ satisfying
$\phi_k=f_k\phi$  for all $k\in J$. The morphisms $f_k$ are called
{\it limit morphisms}. 
(For the example of the category of compact groups see
e.g.\ \cite{compbook},  Definitions 1.25 and 1.27,  or see
Chapter 1 of $\cite{probook}$. For the general concept 
of a limit see \cite{compbook},
Definition A3.41, or go to MacLane's general source book \cite{mac}.)
We have already seen a concrete example of a projective limit in Theorem 
\ref{1.1}. In fact, that example was particular insofar the limit morphisms
$f_k$ were all quotient morphisms. 
 To mathematicians working on the topological algebra of
locally compact groups, projective limits are utterly familiar
by the  Theorem of Yamabe saying that 

\nin
{\it every locally compact topological group $G$ with 
the identity component $G_0$ is a projective limit of
Lie groups provided that  $G/G_0$ is compact.}

\nin (See the classic of 1955 by Montgomery and Zippin \cite{montz}.)

\nin In  particular, this says that every connected locally compact group 
is approximated by connected Lie groups.
Therefore we need to pinpoint in functorial terms
what important theorems like
these say on the principle of ``approximating complicated topological algebraic
structures'' by simpler ones.

\msk

Topologists like to use the concept of a {\it net} on a set $X$ 
generalizing that of
a sequence \cite{kell}: A net $(x_j)_{j\in J}$ is a function $j\mapsto x_j:J\to X$
for a directed poset $J$. 
If $X$ is a topological space and $Y$ a subset of $X$ such that for
every $x\in X$ there is a net $(y_j)_{j\in J}$ of elements in $Y$ such that
$x=\lim_{j\in J}y_j$, then we say that $Y$ is dense in $X$.

\msk
 
So let us now look at a category $\cB$ with a subcategory $\dB$.

\begin{Definition} \label{prolimit}
We call $\dB$ {\it topologically dense} in $\cB$ if it is a full subcategory
of $\cB$ such that for each object $B$ in $\cB$ there is a directed set $J$ and some
projective system 
\vskip-7pt
\cen{
$\{f_{jk}:B_k\to B_j; (j,k)\in J\times J, j\le k\}$}

\nin
of morphisms in $\dB$ such that in $\cB$ the object
$B$ is  (isomorphic to) the projective limit $\lim_{j\in J}B_j$ of
this system with suitable limit morphisms 

\centerline{$B\cong\lim_{k\in J}B_k\mapright{q_j} B_j$,   $j\in J$.}  
\end{Definition}

\vskip-1\baselineskip

As an example we have seen in Theorem \ref{1.1} that the full
subcategory of finite-dimensional unital algebras $\AA_d$ is topologically dense
in the category of weakly complete unital algebras $\AA$.
In the same spirit, by Peter and Weyl,
the category of compact Lie groups is topologically 
dense in the category of compact groups and 
continuous group morphisms (see \cite{compbook},Corollary 2.43).
In \cite{compbook} this Density Theorem is exploited widely.

We owe our readers an explanation of our choice of terminology of
a {\it topologically dense subcategory} which, as we have argued intuitively,
is indeed close to the geometric idea of a dense subspace in a topological space.
The necessity of a comment arises from the fact that in category theoretical
circles, the choice of the terminology of a ``dense subcategory of a category''
is half a century old or older as can be seen from MacLane's standard
text of 1971, where the terminology is introduced close to the
end of the book \cite{mac} on pp.~241, 242, 243. However, that generation of ground
breaking category theoreticians had a distinct leaning towards examples 
supplied by combinatorics and algebra. Therefore, in their eyes,
a category $D$ is, firstly, dense in a category $C$ if every object
of $C$ is a {\sc colimit} of a subsystem of objects from $D$. We would
accordingly suggest to call their approach an approach to {\sc codensity}.
However,  
secondly, their formation of colimits is  {\sc not restricted to
directed systems} (in the way  we insist to use projective limits
 when we (truly!) use limits).  
As a consequence in their terminology, in the category of sets a 
category consisting of one singleton object is codense in the whole category,
and the category consisting of the object $\Z^2$ is codense in 
the category of abelian groups. So dualizing their approach via
Pontryagin would yield that the category consisting of the single object 
of the  traditional torus $\T^2$ would be dense in the category of all compact 
abelian groups. --- At any rate, this predicament causes us to set off our own terminology
of  ``topologically dense subcategories.''

\bsk

\subsection{Density and Adjunction}
\msk

For the class of objects of a category $\cA$ we write $\ob(\cA$).
Now let $L_o\colon \ob(\cA)\to \ob(\cB)$ be a function and $R\colon \cB\to \cA$ 
 a functor and assume that $\cB$ has  a  subcategory $\dB$.

\begin{Definition} \label{basic}
We say that $L_o$ is {\it conditionally left adjoint to $R$
 with respect to a  subcategory $\dB$ of} $\cB$ if for each
$A\in\ob(\cA)$  there is an $\cA$-morphism
$\eta_A\colon A\to RL_o(A)$ such that  for each $B\in \ob(\dB)$ and
each morphism $f\colon A\to R(B)$ in $\cA$
 there is a unique morphism $f'\colon L_o(A)\to B$ in $\cB$
such that $f=R(f')\circ\eta_A$.
\end{Definition} 

A special case illustrates this technical concept:

\begin{Remark} \label{wichtig} If $L_o$ is conditionally left adjoint to $R$
with respect to $\cB$ itself (in place of $\dB$),  
then $L_o$ is the restriction to the objects of
a functor $L\colon \cA\to \cB$
 which is left adjoint to the functor $R$.
\end{Remark}
\begin{Proof} \cite{compbook}, Theorem A3.28.
\end{Proof}

But now we show that the much  weaker condition in Definition \ref{basic}
suffices frequently  for $L_o$ to extend to a left adjoint of $R$.
\bsk
 
\begin{Theorem} \label{density} {\rm (The Density and Adjunction Theorem)}
Assume that  $\cA$ and $\cB$ are two categories and  that $\cB$ has
a topologically dense subcategory $\dB$.
Further assume that  \vskip-10pt
$$
\begin{matrix}
L_o&{\colon}& \ob(\cA)&\to&\ob(\cB)&\mbox{ is a function and}\\
  R&{\colon}&\cB&\to&\cA&\mbox{ is a functor.}\hfill\\
\end{matrix}$$
\vskip-7.5pt
\nin Then the following conditions are equivalent:
\begin{enumerate}[\rm(a)]
\vskip-7.5pt
\item  $L_o$ is conditionally left adjoint to $R$ with respect to $\dB$,
and $R$ preserves projective limits.
\vskip-7.5pt
\item  $L_o$ extends to a left adjoint $L$ of $R$.
\end{enumerate}   
\end{Theorem}

\begin{Proof} For (b) $\Rightarrow$ (a) we refer to 
Remark \ref{wichtig} and to \cite{compbook},
Theorem A3.52, saying that right adjoints are continuous, that is,
preserve all limits. 

Now we prove (a) $\Rightarrow$ (b): In view of Remark \ref{wichtig} it
suffices to show that $L_o$ is conditionally left adjoint to $R$ with
respect to $\cB$ (in place of merely to $\dB$).
So  assume now   that  $A$ and $B$ are objects of
$\cA$ and $\cB$, respectively, and that $f\colon A\to R(B)$ is a morphism
in $\cA$. Then since $\dB$ is topologically dense in $\cB$ we know that 
there exists a projective system

\cen{
$\{f_{jk}\colon B_k\to B_j;\  (j,k)\in J\times J,\  j\le k\}$
of morphisms in $\dB$}

\nin  for some directed set $J$ in $\dB$ such that
$$B=\lim_{j\in J} B_j.\leqno (1)$$ 
Then we obtain a projective system 

\centerline{$\{R(f_{jk})\colon R(B_k)\to R(B_j);\ (j,k)\in J\times J,\ j\le k\}$
of morphisms in $\cA$}

\nin for our directed set $J$. Since $L_o$ is conditionally adjoint to $R$
with respect to $\dB$, for each $j\in J$, then
 there is a unique $\cB$-morphism 
$(Rf_j\circ f)'\colon L_o A\to B_j$
 such that   \def\({\big(} \def\){\big)}
$$Rf_j\circ f=R\((Rf_j\circ f)'\)\circ \eta_A. \leqno(2)$$

\nin
We claim that for $j\le k$ in $J$ we have
$$(Rf_j \circ f)'= f_{jk}\circ (Rf_k\circ f)'.\leqno(3)$$
  For a proof of this claim, we recall from (2) 
that $(Rf_k\circ f)'$
is the {\em unique} $\cB$-morphism for which 
$R\((Rf_k \circ f)'\) \circ \eta_A=Rf_k\circ f$.
Now  
$R\((f_{jk}\circ (Rf_k\circ f)'\)\circ \eta_A$\\
$\begin{matrix}=Rf_{jk} \circ R\((Rf_k\circ f)'\)\circ\eta_A&\mbox{(since $R$ is a functor)}\\
=Rf_{jk}\circ Rf_k\circ f\hfill&\mbox{(by (2))}\hfill\\
=R(f_{jk}\circ f_k)\circ f=Rf_j\circ f&\mbox{(since $R$ is a functor)}\\
=R\(R(f_j\circ f)'\)\circ \eta_A\hfill&\mbox{(by (2)).}\hfill\end{matrix}$

\nin By the uniqueness in the definition of $(Rf_j\circ f)'$
this proves the Claim.

\msk
 By the universal property of the limit, there is now
a unique $\cB$-morphism $f'\colon L_o A\to B$ so that
$$(\forall j\in J)\,(Rf_j\circ f)' =f_j\circ f'. \leqno(4) $$
Consequently, since $R(Rf_j\circ f)'\circ\eta_A=Rf_j\circ f$ by
(2), we have 
$$(\forall j\in J)\, Rf_j\circ f =Rf_j\circ (Rf' \circ\eta_A). \leqno(5) $$
By (a) we  know that we may write $RB=\lim_{j\in J} RB_j$ 
with $Rf_j\colon RB\to RB_j$ as limit morphisms.
By the uniqueness in the universal property of 
the limit (as specified in great generality in
\cite{compbook},  Definition A3.41), from (5) we conclude
$$ (\forall f\colon A\to RB)(\exists! f'\colon L_o A\to B)\quad  
f= Rf'\circ \eta_A.\leqno(6)$$
This completes the proof of (b) \end{Proof}

\bsk

\subsection{An application: The weak completion of a $\K$-vectorspace}

\msk

As an example, consider the functor $W\to \un W\colon \W \to \V$ which assigns to a
weakly complete vector space $W$ the underlying $\K$-vector space
$\un W$. This functor  has a left adjoint $L\colon \V\to \W$ characterized by
the usual universal property recognized in the
usual diagram:

\vskip-15pt 

$$\begin{matrix}& \V&&\hbox to 7mm{} &\W\cr 
\noalign{\vskip3pt}
\noalign{\hrule}\cr
\noalign{\vskip3pt}%
   V&\mapright{\epsilon_V}&\un{L(V)}&\hbox to 7mm{} &L(V)\\
\lmapdown{\forall f}&&\mapdown{\un{L(f')}}&\hbox to 7mm{}&
         \mapdown{\exists! f'}\\
 \un W&\mapright{\id}&\un W&\hbox to 7mm{}&W.
\end{matrix}
$$

\vskip7pt

\nin
The function $f\mapsto f': \V(V,\un W)\to \W(L(V),W)$ is a natural
bijection. 

\bsk

\begin{Proposition} \label{bidual} For a $\K$-vector space $V$ we have
$$L(V)={\un{(V^*)}}^* \mbox{ and }\qquad (\forall \omega\in V^*)\,
 \epsilon_V(v)(\omega)=\omega(v)).$$
\end{Proposition}

\nin Note: In a loose fashion we might write $L(V)=V^{**}$ and say:
\hfill\break 
{\it the weak completion of a $\K$-vector space $V$ is its bidual $V^{**}$}.

\begin{Proof} First we test the universal
property of $L(V)$ for $W\in\ob(\W)$ with $\dim W<\infty$. 
Then the natural morphism $\epsilon_W\colon W\to W^{**}$  
is an isomorphism and for $f\colon V\to W$ we have a commutative diagram
$$\begin{matrix} V&\mapright{\epsilon_V}&V^{**}\\
\lmapdown{f}&&\mapdown{f^{**}}\\
W&\lmapright{\epsilon_W}&W^{**}.\end{matrix}$$
Any $\V$-morphism $f\colon V{\to}W$ yields a
unique morphism  $f'=\epsilon_W^{-1}\circ f^{**}\colon V^{**} \to W$.
The equation $f{=}\un{(f')}\circ\epsilon_V$ is now clear. Thus 
the function $V{\mapsto}{\un{(V^*)}}^*:\ob\V{\to}\ob\W$ has the universal 
property of a 
conditional left adjoint of the functor $W\mapsto\un W$ with respect
to the topologically dense subcategory $\W_d$ of $\W$ consisting of all
 finite-dimensional vector spaces.
So Theorem \ref{density} applies and proves that 
$V\mapsto {\un{(V^*)}}^*$ is left adjoint to $W\mapsto \un W$. 
\end{Proof} 

\rm We note that the necessity of invoking Theorem \ref{density}
indicates that the proof is not entirely trivial.
For $\K=\R$ we have seen that $V^*$ for $V\in\ob(\V)$ is naturally
isomorphic to the Pontryagin dual $\hat V=\Hom_{\rm continuous}(V,\T)$
with $\T=\R/\Z$ when $V$ is endowed with the finest locally convex
topology. (See \cite{compbook}, Theorem A7.10). If, for a 
$W\in\ob(\W)$ we let $W_f$ denote the underlying vector space of $W$
endowed with its finest locally convex topology. Then we have
$$L(V)={\un{(V^*)}}^*\cong (\hat V_f)\hat{\phantom x}.$$
If $A$ is any abelian topological group and ${\hat A}_d$ is its
character group endowed with the discrete topology, then the
compact ``bidual'' $\alpha(A):=(\hat A_d)\hat{\phantom f}$ together
with natural continuous morphism $A\to \alpha(A)$ is the
so called {\it almost periodic compactification} of $A$. We
mention this  here in order to exhibit the analogy between
the weak completion and the almost periodic compactification.

\bsk

\subsection{Strict density and the preservation of projective limits}

\msk

We continue  with categories $\cA$ and $\cB$ having  topologically
dense subcategories
$\dA$, respectively, $\dB$,
 and we consider 
a pair of adjoint functors $L\colon\cA\to \cB$ and $R\colon\cB\to\cA$
between them. Thus we have the following situation

\msk

\nin{\bf Strict Density.} \label{strictdensity} {\it For
each object $A$ of $\cA$ we have some family 
  $q_j\colon A\to A_j$, $j\in Q(A)$ of morphisms with $A_j$ in $\dA$ with
a directed set $Q(A)$ of indices, together with 
a projective system in $\dA$, say,
$q_{jk}\colon A_k\to A_j$ for $j\le k$ in $Q(A)$  such 
that $q_j=q_{jk}\circ q_k$
for $j\le k$, giving us a unique isomorphism $q_A\colon A\to\lim_{j\in Q(A)}A_j$ 
such that
$$\begin{matrix} A&\mapright{q_A}& \lim_{k\in Q(A)}A_k\\
\lmapdown{q_j}&&\mapdown{\rho_j}\hfill\\
A_j&\lmapright{q_{jj}=id}&A_j\hfill\\ \end{matrix}$$
commutes for each $j\in Q(A)$ for the limit morphisms $\rho_j$.

\msk The universal property of the limit will now provide us with
the existence of a crucial morphism
$$ \phi_A\colon L(A)\to \lim_{j\in Q(A)} L(A_j) \leqno(\bullet)$$}
 
\nin We shall investigate this situation in more detail in the remainder of the chapter.

\begin{Lemma} \label{notation}
If $L\colon\cA\to\cB$  is any functor into a complete category $\cB$,
then  
$\{L(q_{jk})\colon L(A_k)\to L(A_j);\ (j,k)\in  Q(A)\times Q(A),\ j\le k\}$
is a projective system in $\cB$, which has a limit 
$$L^\#(A)\defi \lim_{j\in Q(A)} L(A_j)$$
and which provides a morphism $\phi_A\colon L(A)\to L^\#(A)$ 
such that 
$$ \begin{matrix} L(A)&\mapright{\phi_A}&L^\#(A)\\
\lmapdown{Lq_j}&&\mapdown{\rho_j}\\
L(A_j)&\lmapright{=}&L(A_j)\\ \end{matrix}\leqno(7)$$
commutes for all $j\in Q(A)$ for the limit morphisms $\rho_j$.
\end{Lemma}
\bsk

If $\cA$, for example, is the category $\W$ of weakly complete
vector spaces $V$ over $\K=\R$ or $\K=\C$, then a vector subspace $E$ of $V$ is
called a {\it cofinite-dimensional} vector subspace of $V$ if
$\dim V/E<\infty$. Now 
each $V$ defines naturally the filter basis
$J(V)$ of cofinite-dimensional
closed vector subspaces $W$ such that $V\cong\lim_{W\in J(V)} V/W$
 where $V/W$ ranges through the finite-dimensional quotient spaces of $V$
so that the subcategory $\fW$ of all finite-dimensional vector spaces is
topologically dense in $\W$.

\msk

The left adjoint functor $L\colon\cA\to\cB$ preserves
colimits. But it would be  interesting to know whether it preserves also
at least some of the significant limits in the contexts that interest us.
For instance: If $\cA=\W$, the category of weakly complete
$\K$-vector spaces: does then  $L$ under certain 
circumstances preserve projective
limits in $\W$ such as $\lim_{W\in J(V)}V/W$? That is: Is the natural
morphism $\phi_V\colon L(V)\to\lim_{W\in J(V)}L(V/W)$ 
an isomorphism for certain categories $\cB$? 

\bsk
In the example of the category $\W$ of weakly complete vector spaces
each object $V$ gave rise to the projective system 

\centerline{
$\{f_{UW}\colon V/W\to V/U; U,W\in J(V), W\subseteq U\}$,}

\nin
whose limit was naturally isomorphic to $V$.
Here the filter basis $J(V)$ converges to $0$ in the topological space
underlying $V$. If $F$ is any finite-dimensional $\K$-vector space, then
the filter base of vector subspaces $\{F\cap W: W\in J(F)\}$ of $F$
converges to zero in $F$. Now, since $\dim_\K F<\infty$, there is some
member $W_F\in J(V)$ such that $W\subseteq W_F$ implies $F\cap W=\{0\}$.
In terms of quotient morphisms of $V$ this can be expressed as follows:

\msk

Each object $V$ of $\W$ has canonical projective system  of quotient maps
$q_W\colon V\to V/W$, $W\in J(V)$ and $\dim V/W<\infty$
 such that $V\cong \lim_{W\in J(V)}V/W$, and that for each morphism
$f\colon V\to F$ into a finite-dimensional  vector space $F$
we have $\ker f\in J(V)$ so that for every $W\in J(V)$
with $W\subseteq \ker f$ the morphism $f$ factors through
$q_W\colon V\to V/W$. In other words, there 
is an index $W_f\in J(V)$ such that for every
$W\in J(V)$ such that $W\subseteq W_f$ 
there is a morphism
$p_W\colon V/W \to F$ such that $f=p_W\circ q_W$, as in the following
commutative diagram:
$$\begin{matrix}V&\mapright{=}&A\\
\lmapdown{q_W}&&\mapdown{f}\\
V/W&\lmapright{p_W}& F.\end{matrix} $$

\bsk

Following this example we formulate the following definition:

\begin{Definition} \label{strict} For an object $A$ of $\cA$, a projective
system 
$$\{q_{jk}\colon A_k\to A_j: (j,k)\in Q(A)\times Q(A),\quad  j\le k\}\leqno(PS)$$
in $\cA$ will be called {\em appropriate for $A$}
if  $A\cong \lim_{j\in Q(A)}A_j$ with limit morphisms
$q_j\colon A\to A_j$ such that the following
conditions are satisfied at least for a cofinal set of indices $j$ in $Q(A)$:
\begin{enumerate}[\rm(i)]   
\item  For every $\cA$-morphism $f\colon A\to F$ into an $\dA$-object $F$
 there is a $j_0\in Q(A)$ such that for all $j\ge j_0$ there is an $\cA$-morphism 
$p_j\colon A_j\to F$ such that $f=p_j\circ q_j$.

\item  The limit morphisms
$q_j\colon A\to A_j$ are epic. 
\end{enumerate}
If for an object $A$ of $\cA$ there is an 
appropriate projective system $\{q_{jk}\}$, then we say that $A$ is
{\it appropriately representable}.
\end{Definition}
\msk

\begin{Definition} \label{strictlydense}  A subcategory $\dA$ 
of a category $\cA$ is called {\it strictly dense in} $\cA$ 
if  each object $A$ in $\cA$ is  appropriately representable
by a projective system $(PS)$ such that all $A_j$ are
objects from $\dA$.
\end{Definition}

Note that in Condition \ref{strict}(i) the factorisation $f=p_j\circ q_j$
is depicted by the commuting diagram
$$\begin{matrix}A&\mapright{=}&A\\
\lmapdown{q_j}&&\mapdown{f}\\
A_j&\lmapright{p_j}& F.\end{matrix} $$

\ssk

\nin Moreover, condition (ii) is certainly satified if the morphisms
$q_j\colon A \to A_j$ are quotient maps as in the case $\cA=\W$ that we
used as motivation above. Accordingly,
we observe that in the cateory $\W$ of weakly complete
$\K$-vector spaces, the subcategory of finite-dimensional vector spaces 
is strictly dense.

A more sophisticated example is the category of pro-Lie groups, in which
the subcategory of Lie groups is strictly dense.

However, the most relevant example for us is the category $\WA$ 
of weakly complete associative unital algebras in
which the full subcategory $\fWA$ of finite-dimensional $\K$-algebras
is strictly dense by Theorem \ref{1.1}. 

For the applications of the next theorem it is
useful to first recall the following general lemma:

\begin{Lemma} \label{epics} Left adjoint functors preserve epics.
\end{Lemma}
\begin{Proof} By \cite{compbook}, Theorem A3.52 a left adjoint
$L$ preserves colimits. A morphism $e\colon A_1\to A_2$ is an epic if and only if
$$\begin{matrix} A_1&\mapright{e}&A_2\\
    \lmapdown{e}&&\mapdown{\id_{A_2}}\\
             A_2&\lmapright{\id_{A_2}}&A_2\\ \end{matrix}$$

\ssk

\nin is a pushout. A pushout is a colimit (cf. \cite{compbook} EA3.27), 
epimorphisms and  pushouts are dual to  monomorphisms
and  pullbacks; the latter are defined in \ \cite{compbook} Definition A3.9 and
Definition A3.43(ii), respectively.
\end{Proof}

Now we apply Definition \ref{strict}  to provide circumstances in which
the morphism $(\bullet)$ called $\phi_A$ and introduced in 
Lemma \ref{notation} is an epimorphism.

In order to simplify the language of our notation we introduce the 
following definition.

\begin{Definition} \label{suitable pair} A pair of categories $\cA$ and $\cB$ shall be called
a {\it suitable pair of categories} if 
\begin{enumerate}[(i)]
\item $\cA$ posesses a strictly dense subcategory $\dA$,
\item $\cB$ posseses a topologically dense subcategory $\dB$,
\item there is a pair of functors $L\colon\cA \to \cB$ and
$R\colon\cB\to \cA$ such that $L$ is left adjoint to $R$, and
\item $R$ maps $\dB$ into $\dA$.
\end{enumerate}
\end{Definition}

\begin{Theorem} \label{epi} Let $\cA$ and $\cB$ be a suitable 
pair of categories. Assume that the object $A$ of $\cA$ is 
appropriately representable in the form $A=\lim_{j\in Q(A)} A_j$
for an appropriate projective system
$$ \{q_{jk}\colon A_k\to A_j: (j,k)\in Q(A)\times Q(A), j\le k\}.$$

 Then the following statements hold:

$$\{L(q_{jk})\colon L(A_k)\to L(A_j): (j,k)\in Q(A)\times Q(A),\quad  j\le k\}
\leqno({\rm a})$$
is appropriate for $L^\#(A)\defi\lim_{j\in J} L(A_j)$ in $\cB$.

\msk 

\nin{\rm(b)} The morphism 
$$\phi_A\colon L(A)\to L^\#(A)\leqno(\bullet)$$
is an epimorphism.
\end{Theorem}  
\begin{Proof} For proving (b) let $\alpha, \beta\colon L^\#(A)\to B$ be $\cB$-morphisms
such that $\alpha\circ\phi_A=\beta\circ\phi_A$. We must show that 
$\alpha=\beta$. We shall first argue that we may assume that $B$ is in $\dB$.
Since $\dB$ is topologically dense in $\cB$ by \ref{suitable pair}(ii), there is a projective system 
$$\{r_{mn}\colon B_n\to B_n;\ (m,n)\in Q(B)\times Q(B),\ m\le n\}$$
in $\dB$ such that  $B=\lim_{m\in Q(B)} B_m$
with limit morphisms $r_m\colon B\to B_m$ such that 
$r_m=r_{mn}\circ r_n$ for $m\le n$. Then  for each $m\in Q(B)$
we have morphisms 
$r_m\circ\alpha,\ r_m\circ\beta\colon L^\#(A)\to B_m$ such that
$$r_m\circ\alpha\circ\phi_A=r_m\circ\beta\circ\phi_A.\leqno(8)$$
If we can show that for all $m$ we have 
$r_m\circ\alpha=r_m\circ\beta\colon L^\#(A)\to B_m$, then 
by the uniquenes of the universal property of the limit this
will show $\alpha=\beta$, and we shall be done.
So from now on we shall assume that $B$ is in $\dB$.

\msk
 
\nin (a) For a proof of (a) we shall 
have to prove  that the projective system
$$\{L(q_{jk})\colon L(A_k)\to L(A_j): (j,k)\in Q(A)\times Q(A),\quad  j\le k\}$$
with the limit morphisms $\rho_j\colon L^\#(A)\to L(A_j)$ is appropriate,
that is, for each morphism ${\bf f}\colon L^\#(A)\to B$ for an object $B$ in $\dB$,
and all sufficiently large $j$ there will be morphisms 
${\bf p}_j\colon L(A_j)\to B$ such that 
${\bf p}_j\circ\rho_j= {\bf f}:L^\#(A)\to B$ with
the limit morphisms and $\rho_j\colon L^\#(A) \to L(A_j)$.

So let $B$ be a $\dB$-object and ${\bf f}\colon L^\#(A)\to B$ a $\cB$-morphism.
We define $f\defi {\bf f}\circ \phi_A$.
Then $R(f)\colon RL(A)\to R(B)$ is an $\cA$-morphism into an $\dA$ object
$R(B)$ since $R$ maps $\dB$ into $\dA$. By the hypothesis, that $A$ is 
appropriately represented in the form $A=\lim_{j\in Q(A)}A_j$, the morphism
$$A\mapright{\eta_A} RL(A)\mapright{R(f)} R(B)$$
is a morphism from $A$ to an $\dA$ object for which we find a $j_0\in Q(A)$
such that for all $j\in Q(A)$ with $j_0\le j$ there is a $p_j\colon A_j\to R(B)$
such that $R(f)\circ\eta_A = p_j \circ q_j$. By the universal 
property of the left adjoint $L$ there is a unique $\cB$-morphism
$p_j'\colon L(A_j)\to B$ such that $p_j=R(p_j')\circ \eta_{A_j}$.
The following diagram illustrates the situation:

\vskip-15pt

$$\begin{matrix}& \cA&&\hbox to 7mm{} &&\cB\cr
\noalign{\vskip3pt}
\noalign{\hrule}\cr
\noalign{\vskip3pt}%                                                                                       
   A&\mapright{\eta_A}&RL(A)       &\hbox to 7mm{}& L(A) &\mapright{\phi_A}&L^\#(A)\\
\lmapdown{q_j}&&\mapdown{RL(q_j)}  &\hbox to 7mm{}&\mapdown{L(q_j)}&&\mapdown{\rho_j}\\
 A_j&\mapright{\eta_{A_j}}&RL(A_j) &\hbox to 7mm{}&L(A_j)&\mapright{=}     &L(A_j)\\
\lmapdown{p_j}&&\mapdown{R(p_j')}   &\hbox to 7mm{}&\mapdown{p_j'}&  &\mapdown{p_j'}\\
R(B)&\lmapright{=}&R(B)            &\hbox to 7mm{}&B    &\lmapright{=}&B.\\
\noalign{\vskip12pt}
 &&\hskip-40pt R(p_j')\circ RL(q_j)=R(f)  &\hbox to 7mm{}&p_j'\circ L(q_j)=f. &&\\
\end{matrix}$$

\nin We claim that 
$${p_j}'\circ \rho_j={\bf f}:L^\#(A)\to B$$
in the right half of the diagram. For a proof of this claim
we invoke the functoriality of the limit in the following lemma,
which we consider well understood:
\msk
\begin{Lemma}
Let $\xi_j:\lim_{k\in J}X_k \to X_j$ and $\omega_j:\lim_{k\in J}Y_k\to Y_j$ the corresponding
limit cones of two projective limits and assume that there is a compatible
family $\phi_j\colon X_j\to Y_j$ of morphism such that  for all $j\le k$
the diagrams 
$$ \begin{matrix} Y_k&\mapleft{\omega_{jk}}& Y_j\\
\lmapdown{\phi_k}&&\mapdown{\phi_j}\\
X_k&\lmapleft{\xi_{jk}}&X_j\end{matrix}$$
commute. Then there is a unique morphism $\phi\colon\lim_k X_k\to\lim_k Y_k$
such that $\omega_j\circ\phi=\phi_j\circ \xi_j$ for all $j\in J$, i.e.,
the following diagram commutes 
$$\begin{matrix} \lim_kX_k&\mapright{\phi}&\lim_kY_k\\
\lmapdown{\xi_j}&&\mapdown{\omega_j}\\
X_j&\lmapright{\phi_j}&Y_j.\end{matrix}$$ 
\end{Lemma}

\nin Now we apply this lemma to the special case that the $Y_j$ arise as a 
a constant projective diagram with $Y_j=Y$ and $\omega_{jk}=\id_Y$ 
for all $j\le k$ in $J$,
and $\lim_kY_k=Y$ with $\omega_k=\id_Y$ for all $k\in J$. 
Then $\phi\colon \lim_kX_k\to Y$ agrees with $\phi_j\circ \xi_j$
for all $j$, that is
$$\begin{matrix} \lim_kX_k&\mapright{\phi}& Y\\
\mapdown{\xi_j}&&\mapdown{\id_Y}\\
X_j&\lmapright{\phi_j}&Y\end{matrix}\leqno(\dag)$$
commutes for all $j\in J$. This we apply with 
$J=\{j\in Q(A):j_0\le j\}$, $X_j=L(A_j)$, $Y=Y_j=B$, 
$\xi_j=\rho_j\colon L^\#(A){=}\lim_k X_k{\to}X_j{=}L(A_j)$,
$\phi_j{=}p_j'\colon L(A_j){=}X_j{\to}Y=B$, \break %%remove later,if needed 
$\phi={\bf f}\colon L^\#(A)\to B=Y$. 
Then the commuting of  $(\dag)$ yields exactly $p_j'\circ \rho_j={\bf f}$
for $j_0\le j$ as asserted.
So for all $j\in Q(A)$, $j_0\le j$, the morphisms $p_j'\colon L(A_j)\to B$ are the required 
morphisms ${\bf p}_j\colon L(A_j)\to B$.
\vskip7pt

\msk

(b) We finally   prove that $\phi_A$ is an epic.
By (b) there is a $j_0$ such that for all $j\ge j_0$ there exist 
$\cB$-morphisms $\alpha_j\colon L(A_j)\to B$ and $\beta_j\colon L(A_j)\to B$
such that $\alpha =\alpha_j\rho_j$ and $\beta=\beta_j\rho_j$:
$$\begin{matrix}L^\#(A)&\mapright{\rho_j}&L(A_j)\\
\lmapdown{\alpha}\mapdown{\beta}&&\lmapdown{\alpha_j}\mapdown{\beta_j}\\
B&\lmapright{\id_B}&B.\end{matrix}\leqno(9)$$
Then we must show that 
$$(\forall j\ge j_0)\, \alpha_j=\beta_j.\leqno(*)$$

Now  by the Definition of $\phi_A$ we have 
$$(\forall j\in J)\,\rho_j\phi_A=L(q_j).\leqno(10)$$
We consider the following diagram
$$\begin{matrix}L(A)&\mapright{\id_{L(A)}}&L(A)\\
\lmapdown{\phi_A}&&\mapdown{L(q_j)}\\
L^\#(A)&\mapright{\rho_j}&L(A_j)\\
\lmapdown{\alpha}\mapdown{\beta}&&\lmapdown{\alpha_j}\mapdown{\beta_j}\\
B&\lmapright{\id_B}&B.\end{matrix}\leqno(11)$$

\nin The top square commutes by (10) for all $j\in J$. The two bottom 
squares commute for each $\alpha$ and $\beta$ and all $j\ge j_0$
by (9).  Accordingly, the outside rectangles commute for
both $\alpha$ and $\beta$. The left vertical edges $\alpha\phi_A=\beta\phi_A$
agree by assumption on $\alpha$ and $\beta$. 
So for each $j\ge j_0$ we compute
$$\alpha_jL(q_j)=\id_B\alpha\phi_A\id_{L(A)}^{-1}
=\id_B\beta\phi_a\id_{L(A)}^{-1}=\beta_jL(q_j).\leqno(12)$$
The morphisms $q_j\colon A\to A_j$ are epic by Definition 7.8(ii).
Then Lemma 7.10 shows that the morphisms $L(q_j)\colon L(A)\to L(A_j)$
are all epic. Now (11) implies that $\alpha_j=\beta_j$ for 
all $j\ge j_0$. So $(*)$ is proved and this is what we
had to show. 
\end{Proof}

Notice that Theorem \ref{epi} does not assert that the objects $L(A_j)$
are  (even cofinally) objects of the topologically dense subcategory $\dB$. In fact,
in the applications, which we aim for, this is not the case. It is nevertheless
assumed by Definition \ref{suitable pair}(ii) that every object $B$ of $\cB$ is a projective limit
of objects from $\dB$. 

\begin{Lemma}\label{$^*$}  
Let $\cA$ and $\cB$ be a suitable pair of categories.  
Assume that the object $A$ of $\cA$ is appropriately representable 
as $A=\lim_{j\in Q(A)} A_j$.
Abbreviating
$\lim_{j\in Q(A)}L(A_j)$ by $L^\#(A)$,  define
an $\cA$-morphism $\eta_A^\#\colon A\to R(L^\#(A))$
by $\eta_A^\#=R(\phi_A)\circ \eta_A$.

Then for each object $B\in\dB$ and each 
$\cA$-morphism $f\colon A\to RB$ there is a unique $\cB$-morphism
$f^\#\colon L^\# (A)\to B$ such that $f=R(f^\#)\circ \eta^\#_A$.
\end{Lemma}

 \begin{Proof}  By Definition \ref{suitable pair}(iv), the functor $R$ maps $\dB$
into $\dA$. So $RB$ is in
$\dA$.   By Definition \ref{suitable pair}(i), the subcategory  $\dA$ is strictly
dense in $\cA$. Since $A$ is appropriately representable, there indeed exists  a projective system  
$$\{q_{jk}\colon A_k\to A_j: (j,k)\in Q(A)\times Q(a),\quad  j\le k\}$$
which is  appropriate for $A$. 
 So  there is a $j_0\in Q(A)$ such that for all $j$ with $j_0\le j$ 
there are $\dA$ morphisms
$p_j\colon A_j\to R(B)$ such that $f=p_j\circ q_j$. Since 
$L$ is left adjoint to $R$, there are unique $\cB$-morphisms
$f'\colon LA\to B$ and $(p_j)'\colon LA_j\to B$ such that
$f=Rf' \circ \eta_A$ and $p_j=R(p_j)'\circ\eta_A$. Now
from $f=p_j\circ q_j$, by  \cite{compbook}, Proposition A3.33
we deduce
$$f'=(p_j)'\circ Lq_j.\leqno(13)$$
The fill-in
morphism $\phi_A\colon LA\to L^\# A$ of Lemma \ref{notation}
satisfies $Lq_j= \rho_j\circ \phi_A$ with the limit morphism
$\rho_j\colon L^\# A\to LA_j$. We set $f^\#_j=(p_j)'\circ\rho_j$.
If $k\in Q(A)$ satisfies $j\le k$, then we have a commutative 
diagram
$$\begin{matrix} L^\# A&\lmapright{=}&L^\# A\\
\lmapdown{\rho_k}&&\mapdown{\rho_j}\\
L(A_k)&\mapright{L(q_{jk})}&L(A_j)\\
\lmapdown{(p_k)'}&&\mapdown{(p_j)'}\\
B&\mapright{{=}}& B\\ \end{matrix}$$
with $f^\#_j=(p_j)'\circ \rho_j =(p_k)'\circ \rho_k=f^\#_k$.
That is, for a cofinal subset $Q_c(A)\subseteq Q(A)$ of
$Q(A)$ the function $k\mapsto f^\#_k\colon Q_c(A)\to \cB(L^\#(A),B)$
is constant. Hence we have a unique morphism  $f^\#\colon L^\# A\to B$
such that $f^\#=f^\#_j$ for all sufficiently large $j$ such that
 (7) and (13) imply  
$$f^\#\circ \phi_A=(p_j)'\circ\rho_j\circ\phi_ A=(p_j)'\circ Lq_j=f'\leqno(14)$$
for all sufficiently large $j$, and thus
$$R(f^\#)\circ\eta_A^\#=R(f^\#)\circ R(\phi_A)\circ\eta_A=R(f')\circ\eta_A=f.$$
If $f^*\colon L^\# A\to B$ is a $\cB$-morphism such that 
$f=R(f^*)\circ\eta_A^\# =R(f^*)\circ R(\phi_A)\circ \eta_A$, then
$f^*\circ\phi_A =f'$ by the uniqueness in determining $f'$.
Since also $f^\#\circ\phi_A=f'$ by (14) above, we may
conclude $f^*=f^\#$, since $\phi_A$ is an epimorphism
by Theorem \ref{epi}.
This completes the proof of Lemma \ref{$^*$}.
\end{Proof}
\msk
As a corollary of the epimorphism Theorem \ref{epi} we
now have the following main result, in which it happens that a left
adjoint functor $L$ preserves, in addition to all colimits, also certain
limits. In its formulation  we retain the notation
of Lemma \ref{notation} and Definition \ref{strict}. 

\begin{Theorem} \label{prolim}
Let $\cA$ and $\cB$ be a suitable pair of categories
and assume that
  the  projective system  
$$\{q_{jk}\colon A_k\to A_j: (j,k)\in Q(A)\times Q(A),\quad  j\le k\}$$
is appropriate for $A$.
Then the morphism
$$\phi_A\colon L(\lim_{j\in Q(A)} A_j)\to \lim_{j\in Q(A)} L(A_j)\leqno{(\bullet)}$$
is an isomorphism.
\end{Theorem}  
 
\begin{proof}  From Lemma \ref{$^*$} and Theorem \ref{density} it now follows
that $L^\#$ extends to a functor $L^\#\colon\cA\to\cB$ which is left adjoint
to $R$. Thus $L$ and $L^\#$ are naturally isomorphic functors.
Then there is a commutative diagram of  natural functions
$$\begin{matrix}
\cB(L^\#(A),B)&\lmapright{\beta_{AB}}&\cB(L(A),B)\\
\lmapdown{\alpha^\#_{AB}}&&\mapdown{\alpha_{AB}}\\
\cA(A,R(B))&\mapright{=}&\cA(A,R(B)),\\ \end{matrix}\leqno(15)$$
\begin{enumerate}[\rm(a)]
\item $\beta_{AB}(h)= h\circ\phi_A$ for $h\colon L^\#(A)\to B$,
\item $\alpha^\#_{AB}(h)=R(h\circ\phi_A)\circ\eta_A$, for $h\colon LA\to B$,
\item $\alpha_{AB}(h)=R(h)\circ\eta_A$, for $h\colon L(A)\to B$.
\end{enumerate}
The bijectivity of $\alpha_{AB}$ expresses the fact that
$L$ is left adjoint to $R$, and likewise the bijectivity of
$\alpha^\#_{AB}$ is now secured since we proved that
$L^\#$ is left adjoint to $R$. the commutativity of the diagram
(15) then shows the bijectivity of $\beta_{AB}$ which in turn
proves that $\phi_A$ is an isomorphism. This completes the
proof.
\end{proof}

Since the right adjoint $R\colon \cB\to \cA$ preserve limits, the
following corollary is immediate:

\begin{Corollary} \label{prolim2}  \hskip-8pt Under the hypotheses of  
\hskip5pt{\rm Theorem \ref{prolim}}, for each $A{=}\lim_{j\in\Q(A)}A_j$,
 the $\cA$-morphism
$R(\phi_A)\colon RL(A)=RL(\lim_{j\in Q(A)} A_j)\to \lim_{j\in Q(A)} RL(A_j)$
is an isomorphism.
\end{Corollary}

\nin
If $r_j\colon L^\#(A)\defi\lim_k L(A_k)\to L(A_j)$, $j\in Q(A)$ denotes
the limit morphisms, the situation is illustrated by the following
diagram: 

\vskip-1\baselineskip
$$
\begin{matrix}&&\cA&&&\hbox to 25mm{} &\cB\cr 
\noalign{\vskip3pt}
\noalign{\hrule}\cr
\noalign{\vskip3pt}%
A&\mapright{\eta_A}&RL(A)&\mapright{R(\phi_A)}&\lim_k RL(A_k)&\hbox to 25mm{} &L(A)\\
\lmapdown{\forall q_j}&&\mapdown{RL(q_j')}&&\mapdown{R(r_j)}&\hbox to 25mm{}&
         \mapdown{\exists! q_j'}&&\\
A_j&\lmapright{=}&A_j&\lmapright{\eta_{A_j}}&RL(A_j)&\hbox to 25mm{}&L(A_j)\\
\end{matrix}
$$  

\medskip

\begin{Corollary} \label{prolim3} Assume the hypotheses of   
{\rm Theorem \ref{prolim}}, and, in addition, that for all objects
$A\in\ob(\cA_d)$ the front adjunction $\eta_A\colon A\to RL(A)$ is monic.
Then it is monic for all objects $A\in\ob(\cA)$ in $\cA$.
\end{Corollary}
\begin{Proof} Let $\alpha, \beta\colon X\to A$ be morphisms such that
$\eta_A\alpha=\eta_A\beta$. Then for $j\in Q(A)$ we have
$\eta_{A_j}q_j\alpha=\eta_{A_j}RL(q_j')\eta_A\alpha
=\eta_{A_j}RL(q_j')\eta_A\beta=\eta_{A_j}q_j\beta$.
Since $\eta_{A_j}$ is monic, we have
$$(\forall j\in Q(A))\, q_j\alpha =q_j\beta.$$
Since $A=\lim_{j\in Q(A)} A_j$, the uniqueness of the morphism in the universal
property of the limit (see. e.g.\ \cite{compbook}, Definition  A3.41) 
implies $\alpha=\beta$.
\end{Proof} 
 
\bsk 

\subsection{Application: $\cB=\AA$.}

\msk

Our main target category for various left adjoint functors is the
category $\AA$ of weakly complete unital algebras over $\K=\R,\C$.
\begin{Proposition}\label{anwendung} The subcategory $\AA_d$ of finite-dimensional
associative unital $\K$-algebras is strictly dense in $\AA$.
\end{Proposition}
\begin{Proof} Each weakly complete algebra $A$ has a filter base $Q(A)$ of
closed two sided ideals $I$ such tht $A/I$ is a finite-dimensional $\K$ algebra,
and the natural morphims $q_A\colon A\to \lim_{I\in Q(A)} A/I$ is an
isomorphism according to Theorem \ref{1.1}. Thus by Definition \ref{prolimit},
$\AA_d$ is topologically dense in $\AA$.

We need to verify the conditions of Definition \ref{strict}. Let $f\colon A\to F$
be an $\AA$-morphism for a finite-dimensional algebra $F$ and let $I=\ker f$.
Then $A/I\cong\im f$  (see \cite{compbook}, Theorem  A7.12(b)) and so,
since $\dim \im f\le \dim F<\infty$,  $I\in Q(A)$. Let $q_I\colon A\to A/I$ denote
the quotient morphism and $p_I\colon A/I\to F$ the injective morphism
induced by $f$. Then $f=p_I\circ q_I$. Hence condition (i) of 
Definition \ref{strict} is satisfied.
Since the quotient morphisms $q_I\colon A\to A/I$ are surjective and therefore
epimorphisms,  condition (ii) is satisfied as well. 
This completes the proof.  \end{Proof}

For observing first applications, we let $\cA$  be a category $\TA$ of topological 
algebraic structures  and $R\colon \AA\to \TA$ a limit preserving functor.
We assume that $R$ satisfies the Solution Set Condition
(see \cite{compbook} A3.58), as is the case in the examples we discuss below
(cf. \cite{dhtwo,hofkra}). 
Hence a left adjoint functor $L\colon \TA\to \AA$ exists (see
\cite{compbook}, Theorem A3.60.)  

\msk

Assume that $L\colon\TA\to
\AA$ is left adjoint to $R$ and
that $\eta_X\colon X\to RL(X)$ is the front adjunction 
(see \cite{compbook}, Definition A3.37). 
Now $L(X)$ is a weakly complete unital algebra. In practically all
examples of interest to us, for a weakly complete unital algebra $A$,
the $\TA$-object $R(A)$ is a subset of $A$ such as for instance $A^\times$ (if
$\TA$ is a category of topological groups), or $\li A$ (if $\TA$ is a
category of topological Lie algebras), or the underlying topological
 space $|A|$ (if $\TA$ is a category of topological spaces).

 In such a situation $\eta_X\colon X\to RL(X)$ is a function and we can
consider its image $\eta_X(X)$ as a subset of  $RL(X)$. Then 
$\<\eta_X(X)\>$ denotes  the smallest unital
subalgebra  containing $\eta_X(X)$ and $\overline{\<\eta_X(X)\>}$
the smallest $\AA$ subobject of $L(X)$.
Under these circumstances we define a function 
$$L_o\colon\ob(\TA)\to \ob(\AA)\mbox{ by }L_o(X)=\overline{\<\eta_X(X)\>},$$ 
the smallest $\AA$-subobject  for which the morphism
$\eta_X\colon X\to LR(X)$ factors through the inclusion morphism
$R(L_o(X))\to RL(X)$. Then  one observes immediately that
$L_o$ is conditionally left adjoint to $R$
with respect to $\AA$. (Cf.\ Definition \ref{basic}.) However,
here Remark \ref{wichtig} applies and shows that the 
containment $L_o(X)\subseteq L(X)$ is equality in all cases. Thus we
have

\begin{Proposition} \label{epimorphs} Assume that $R\colon\AA\to\TA$
and $L\colon\TA\to\AA$  is a pair of adjoint functors where
$\TA$ is a category of topological algebraic structures for which the front adjunctions
$\eta_X\colon X\to RL(X)$ are functions whose image $\eta_X(X)$ is
a subset of the weakly complete unital algebra $L(X)$. Then
for each $X\in\ob(\TA)$,
the abstract unital algebra $\<\eta_X(X)\>$ generated by the 
image of $\eta_X$ is dense in the weakly complete unital algebra $L(X)$.
\end{Proposition}

\msk

Our immediate 
examples for the category $\TA$ are as follows:

\msk

\nin{\bf(A)}\  $\TA = \progr$, the category of pro-Lie  groups, $RA=A^\times$
the group of units of the weakly complete algebra $A$. For the
fact that $A^\times$ is a pro-Lie group see \cite{dhtwo} or \cite{compbook},
Proposition A7.37. The left adjoint L is the 
{\it weakly complete group algebra  $G\mapsto \K[G]$ over $\K$}. 
It was discussed in \cite{dhtwo}, \cite{hofkra}, and \cite{compbook} (mostly
for $\K=\R$ and compact groups $G$). 

\ssk

A prominent subcategory of $\progr$ is the full subcategory
$\compgr$  of compact groups for which the real weakly complete
group algebra $\R[G]$ is particularly effective. See \cite{dhtwo}.

\msk

\nin{\bf(B)}\  $\TA = \prolie$, the category of 
profinite-dimensional Lie algebras
over $\K$ (cf.\ \cite{hofkra}, \cite{hofkraii}). The functor 
$R\colon \AA \to\prolie$ associates with a weakly complete unital algebra
$A$ the profinite-dimensional Lie algebra $\li A$ defined on
the weakly complete underlying weakly complete $\K$-vector space
endowed with the Lie algebra multiplication $[x,y]=xy-yx$.
 Then the left adjoint $L\colon \prolie \to\AA$ is the
{\it weakly complete universal enveloping algebra $\g\mapsto \UU_\K(\g)$
over  $\K$} \cite{hofkra,hofkraii} which we shall address again below.

\msk

\nin{\bf(C)}\ $\TA= \W$, the category of weakly complete $\K$-vector spaces.
The functor $R\colon \AA\to \W$ associates with a weakly complete
unital algebra $A$ the underlying weakly complete topological 
$\K$-vector space $|A|$.
The left adjoint $L\colon \W\to \AA$ of $R$ is, as we shall discuss
in the subsequent section, the functor which associates with
any weakly complete vector space $W$    the {\it weakly complete
tensor algebra $\TT(W)$ of $W$ over $\K$}.

\msk

\begin{Proposition} \label{morphs} {\rm(i)} In each of the categories 
$\TA= \compgr$, $\prolie$, and $\W$, a monomorphism $f\colon X\to Y$
induces an isomorphism $X\to f(X)$ onto the image, that is, an embedding
in the respective category $\TA$. 

{\rm(ii)} The front
adjunction $\eta_X\colon X\to RL(X)$, namely,

\cen{$G\to \K[G]^\times$, $\g\to \li{\UU_K(\g)}$, and $W\to |\TT(W)|$,}

\nin  is an embedding in the respective category. 
That is, $X$ may be considered as a $\TA$-subobject of $RL(X)$
and a subset of $L(X)$. 

{\rm(iii)} If $X\in\ob(\TA)$ and $X\subset L(X)$ as in {\rm(ii)}
above, then $\<X\>$, the abstract unital algebra generated by
$X$ in $L(X)$ is dense in $L(X)$.
\end{Proposition}

\begin{Proof} Part (i) may be safely considered as an exercise;
for the two categories $\prolie$ and $\W$ see also \cite{compbook}, 
Theorem A7.12.

Part (ii) is then a consequence of Part (i), 
Corollary \ref{prolim3} and the following 
facts which secure that the front adjunction 
$\eta_X\colon X\to RL(X)$ is injective,
hence monic for $X\in\ob(\TA)$:

Part(iii) follows from Proposition \ref{epimorphs}.
\ssk

(a) Every compact Lie group has a faithful linear representation (see e.g.
\cite{compbook}, Corollary 2.40).
 
\ssk

(b) Every finite-dimensional Lie algebra over a field of characteristic 0
has a faithful linear representation (Ado's Theorem, 
see \cite{bour}, Chap. 1, Paragraph  7, n$^{\rm o}$ 3,  Th\'eor\`eme 3). 

\ssk

(c) It suffices to observe that the one-dimensional vector space $\K$
has a faithful linear representation, e.g.
$$ c\mapsto\begin{pmatrix} 1&c\\ 0&1\end{pmatrix}.$$
\vskip-22pt
\end{Proof}

We now secure the validity of the hypotheses of Theorem \ref{prolim} 
for the examples {\bf (A), (B)}, and {\bf(C)}.

\begin{Proposition} \label{linus-a} For the functors $L=\K[-]$, ${\bf U}_\K$,
and $\bf T$ the morphism $\phi_A\colon A\to\lim_{j\in Q(A)}L(A_j)$
is an isomorphism.
\end{Proposition}

\begin{Proof} We show that the hypothesis of Theorem \ref{prolim} is
satisfied in each of the three examples {\bf (A), (B)}, and {\bf(C)}.

\msk

{\bf(A)} $G\mapsto \K[G]: \progr\to \AA$, is left adjoint to 
$A\mapsto A^\times: \AA\to \progr$. We note that the subcategory
${\cal LIE}$ of Lie groups is strictly dense in $\progr$: 

 Each pro-Lie group $G$
is the projective limit of its Lie group quotients $G/N$, $N\in {\cal N}(G)$,
where ${\cal N}(G)$ denotes the filter basis of normal subgroups $N$
of $G$ such that $G/N$ is a Lie group. 
If $f\colon G\to L$ is a morphism of $G$ into a Lie group, let $N$ be the
kernel of $f$. Then $f$ factors through the quotient morphism $q\colon G\to G/N$
followed by an injection of Lie groups $G/N\to L$. Each quotient morphism
$G\to G/N$ is an epimorphism. So 
$$\{q_{MN}:G/N\to G/M: (M,N)\in{\cal N}(G)\times{\cal N}(G), N\subseteq M\}$$
is appropriate for $G$. The functor $A\mapsto A^\times: \AA\to\progr$ maps 
finite-dimensional algebras to Lie groups. So the hypothesis of Theorem
\ref{prolim} is satisfied and so 
$$\phi_G\colon G\to \lim_{N\in {\cal N}(G)}\K[G/N]^\times$$ 
is an isomorphism.
\msk 
The cases {\bf(B)} and {\bf(C)} are equally simple and are left as an
exercise.
\end{Proof}

\nin Theorem \ref{prolim}  now has immediately the following corollaries:

\begin{Theorem} \label{solution-a} Each pro-Lie group $G$ has an appropriate projective
limit representation $G=\lim_{j\in J} G_j$ in terms of Lie groups. Therefore
$$\K[G]\cong \lim_{j\in J} \K[G_j].$$  
\end{Theorem}

\begin{Theorem} \label{solution-b} Each profinite-dimensional 
weakly complete Lie algebra has an appropriate  projective limit 
representation $\lim_{j\in J} \g_j$ in terms of finite-dimensional
Lie algebras. Therefore 
$$\UU_\K(\g)\cong\lim_{j\in J}\UU_\K(\g_j).$$
\end{Theorem}

\begin{Theorem} \label{solution-c} If a weakly complete $\K$-vector space $W$
is represented in terms of  an appropriate  projective limit representation 
$W=\lim_{j\in J} W_j$ in terms of
finite-dimensional vector spaces. Therefore 
$$ \TT(W)\cong \lim_{j\in J} \TT(W_j).$$
\end{Theorem}

\section{Appendix: The Definition of the Tensor Algebra}

\rm

In Paragraph {\bf(C)} above we already introduced the tensor algebra 
of a weakly complete vector space. Let us now review this concept
more systematically. 
So we let $\K$ again denote one of the topological fields $\R$ or $\C$,
and 
 $\AA$  the category of weakly complete associative 
unital algebras over $\K$.

\msk

Here is the  definition of the tensor algebra via its universal property:

\begin{Theorem}  \label{1.3} {\rm (The Existence Theorem of $\TT$)} 
The underlying weakly complete vector space functor $A\mapsto|A|$
from $\AA$ to $\W$ has a left adjoint $\TT\colon\W\to\AA$.

The front adjunction $\omega_V\colon V\to |\TT(V)|$ is an embedding
of topological vector spaces.\end{Theorem}

\begin{Proof} The category $\W$ is complete. 
(Exercise. Cf.\ Theorem A3.48 of 
\cite{compbook}, p.\ 819.) The
 ``Solution Set Condition'' (of Definition A3.58 in \cite{compbook}, 
p.\ 824) holds.
(Exercise: Cf.\ the proof Lemma 3.58 of \cite{compbook}.) 
Hence $\TT$ exists by the Adjoint Functor Existence Theorem
(i.e., Theorem A3.60 of \cite{compbook}, p.\ 825).

The assertion about $\omega_V$ follows from 
Proposition \ref{morphs} (iii).\end{Proof}

In other words, 
 each weakly complete vector space $V$ may be considered as a weakly 
complete vector subspace of the weakly complete tensor algebra
$\TT(V)$ with the property  that  each continuous 
linar map $f\colon V\to|A|$ with some weakly complete associative unital 
algebra $A$ and its underlying weakly complete vector space $|A|$
extends uniquely to a $\AA$-morphism $f'\colon \TT(V)\to A$. 

$$\begin{matrix}& \W&&\hbox to 7mm{} &\AA\cr 
\noalign{\vskip3pt}
\noalign{\hrule}\cr
\noalign{\vskip3pt}%
   V&\mapright{\subseteq}&|\TT(V)|&\hbox to 7mm{} &\TT(V)\\
\lmapdown{\forall f}&&\mapdown{|f'|}&\hbox to 7mm{}&
         \mapdown{\exists! f'}\\
 |A|&\lmapright{\id}&|A|&\hbox to 7mm{}&A.
\end{matrix}
$$

\medskip

\noindent If necessary we shall write $\TT_\K$ instead of $\TT$ 
whenever the ground field should be emphasized.

\begin{Definition} \label{1.4} For each weakly complete $\K$-vector space $V$
 we shall call $\TT_\K(V)$ 
{\it the weakly complete tensor algebra} of $V$ 
(over $\K$).
\end{Definition}

We record what we already saw in Section 1 in %Proposition \ref{morphs} and
Theorem \ref{solution-c}:

\begin{Corollary} \label{tensor}  If $V$ is represented as a 
projective limit $\lim_{j\in J} V_j$
of finite-dimensional vector spaces, then $\TT(V)\cong\lim_{j\in J}\TT(V_j)$.
\end{Corollary}

Every unital associative algebra $A$ has injective
morphism $\iota_A\colon \K\to A$ given by $\iota_A(t)=t\.1$. 
In some circumstances, $\iota$ is a coretraction:

\msk

\begin{Remark} \label{Remark-2.X}  For every weakly complete vector space $V$, the morphism \break
$v\mapsto 0: V\to \K$, according to the definition of $\T(V)$, induces
a natural $\W$-morphism $\alpha_V\colon {\bf T}(V)\to\K$ such that 
$\alpha_V(V)=\{0\}$ and that $\alpha_V\circ\iota_{{\bf T}(V)}=\id_\K$.     
\end{Remark}
\msk

The retraction  $\alpha_V$ is frequently called the {\it augmentation}.
of ${\bf T}(V)$.

\medskip Let us compare the weakly complete tensor 
algebra with the abstract tensor
algebra $T(E)$ of a plain $\K$-vector space $E$.  
By the universal property of the abstract tensor product
$T(|V|)$, the linear inclusion map $\xi_{|V|}\colon |V|\to T(|V|)$ extends
to a unique morphism of unital algebras $j_{|V|}\colon T(|V|)\to |\TT(V)|$
such that 
$$\begin{matrix} |V|&\mapright{\xi_{|V|}}&T(|V|)\\
\lmapdown{|\omega_V|}&&\mapdown{j_V}\\
|\TT(V)|&\lmapright{\id}&|\TT(V)|\\\end{matrix}\leqno(1)$$
commutes.

\medskip For a natural number $m$ and a weakly complete vector space $V$,
set $A_m=\bigotimes^m_\W V$. Then $A_m$ is a weakly complete vector space.
If $m$, and $n$ are natural numbers, then 
there is a canonical continuous bilinear
map
$$(a_m,a_n)\mapsto a_ma_n:=a_m\otimes_\W  a_n:A_m\times A_n\to A_{m+n},$$  
where we have identified the naturally isomorphic weakly complete
vector spaces $A_m\otimes_\W A_n$ and $A_{m+n}$. Set $A_0=\K$. Then
$\bigoplus_{m=0}^\infty A_m$ is a graded unital algebra $D$ with the 
multiplication 

\cen{$(a_m)_{m\in\N_0}(b_m)_{m\in\N_0}=(\sum_{j+k=m}a_jb_k)_{m\in\N_0}$,}

\nin  dense in
the weakly complete vector space $A(V):=\prod_{m=0}^\infty A_m$.
By the definition of the unital algebra $D$, there is a unique injective
morphism of unital algebrs 
$i_V\colon T(|V|)=\bigoplus_{m=0}^\infty\bigotimes^m_\V V\to|A(V)|$
so that we have a commutative diagram:
$$\begin{matrix} |V|&\mapright{\xi_{|V|}}&T(|V|)\\
\lmapdown{\id}&&\mapdown{i_V}\\
V&\lmapright{\incl}&A(V).\\\end{matrix}\leqno(2)$$
\ssk 
Now,  multiplication in $D$ is continuous w.r.t.\ the topology induced from
$A(V)$ and therefore extends continuously to a multiplication on $A(V)$,
making $A(V)$ a weakly complete unital algebra. There is an injective
continuous linear map $\iota_V\colon V\to A(V)$ given by
$\iota_V(v)=(0,v,0,0,\dots)\in \K\times A_1\times A_2\times\cdots=A(V)$
which by Theorem \ref{1.3} yields a unique morphism of weakly complete unital algebras
${\iota_V}'\colon T(V)\to A(V)$ such that $\iota_V(v)={\iota_V}'(\omega_V(v))$
for all $v\in V$, i.e.\ such that the following diagram commutes:
$$
\begin{matrix}
   V&\mapright{\omega_V}&|\TT(V)|&\hbox to 7mm{} &\TT(V)\\
\lmapdown{\iota_V}&&\mapdown{|{\iota_V}'|}&\hbox to 7mm{}&
         \mapdown{{\iota_V}'}\\
 |A(V)|&\lmapright{\id}&|A(V)|&\hbox to 7mm{}&A(V).
\end{matrix}\leqno(3)
$$

We now have $j_V\circ\xi_{|V|}=|\omega_V|$ by (1),
$i_V\circ \xi_{|V|}=|\iota_V|$ by (2), and $|\iota_V|=|{\iota_V}'|\circ|\omega_V|$
by (3). Therefore $i_V\circ\xi_{|V|}=|{\iota_V}'|\circ j_V\circ \xi_{|V|}$,
and so the uniqueness in the universal property of $T(|V|)$ allows us
to conclude 
$$i_V=|{\iota_V}'|\circ j_V.\leqno(4)$$
But $i_V$ is injective, and so $j_V\colon T(|V|)\to |\TT(V)|$
is injective.

\msk

Collecting the information we have collected we now
arrive at the following insight:

\begin{Lemma} \label{test} For each weakly complete vector space $V$,
the weakly complete unital algebra $\TT(V)$ contains a  copy of
the algebraic tensor algebra $T(|V|)=\bigoplus_{m=0}^\infty \bigotimes^m|V|$
algebraically generated by $V\subseteq\TT(V)$.
\end{Lemma}
\begin{Proof}  The completion of the proof is now
 an exercise. \end{Proof}
\msk

\begin{Theorem} \label{density2} {\rm(i)} For any weakly complete vector space $V$,
the unital associative subalgebra $\<V\>$ generated algebraically
in $\TT(V)$ by $V$  is dense in $\TT(V)$. 

{\rm(ii)} Moreover, $\<V\>$ is algebraically isomorphic to the algebraic
tensor algebra $T(|V|)$ generated by $|V|$.
\end{Theorem}

\begin{Proof} (i) The assertion was proved in Proposition \ref{morphs}(iii).

(ii) By the universal property of
the algebraic tensor algebra $T(|V|)$ generated by $|V|$ there is a
 morphism $j_V\colon T(|V|)\to |\TT(V)|$ (see (1) above)
whose corestriction to its image is a morphism of unital algebras 
from $T(|V|)$ to $\<V\>$ which is
is injective by Lemma \ref{test} and therefore is an isomorphism of unital
algebras. 
\end{Proof}

\bsk

Let us now use the weakly complete tensor algebra to construct $\UU_\K(\g)$
as a quotient of $T(|\g|)$.
In the classical theory of universal enveloping algebras of Lie algebras,
the construction usually does  proceed from the tensor algebra as an origin and
progresses to the enveloping algebra as a quotient. In the world $\W$ of
weakly complete vector spaces we proceeded systematically via universal
properties using category theoretical standard methods. In this fashion
we have developed the weakly complete tensor algebra and the weakly complete
universal enveloping algebra separately albeit with unified methods. 
Now let us pause
and bring the two together again using the principle of the universal property.

\bsk 

Let $\g$ be a profinite-dimensional Lie algebra and $|\g|$ the weakly complete
vector space on which it is based. Let $\|\g\|$ denote the underlying 
vector space
and $\u{\g}$ the underlying abstract Lie algebra. 
There is a quotient morphism of unital algebras
$p_{|\g|}\colon T(\|\g\|)\to U(\u{\g})$ well known form the apparatus
of the 
Poincar\'e-Birkhoff-Witt-Theorem where we  may consider $\|\g\|$ as a vector
subspace of $T(\|\g\|)$.  Now we elevate this quotient to the level
of the weakly complete unital algebras. From Theorem \ref{1.3a} we know that
$\g\subseteq \UU(\g)$. This give us an embedding of weakly complete
vector spaces  $|\g|\to|\UU(\g)|$ where $\UU(\g)$ is a weakly complete unital
algebra. Then Theorem \ref{1.3} provides us with a unique morphism 
$q_\g\colon \TT(|\g|)\to \UU(\g)$ extending the identity function 
$|\g|\to \g$. As a morphism of weakly complete algebras, $q_\g$ has a closed image
(see e.g.\  \cite{compbook}, Theorem A7.12), and by Theorem \ref{1.3} (i) has
a dense image. Thus $q_\g$ is surjective and thus a quotient map (againby
\cite{compbook}, A7.12). We summarize this in the following Theorem whose proof is clear
from what we know:

\begin{Theorem}\label{completion}
There is a canonical quotient morphism of weakly complete algebras
$q_\g\colon \TT(|\g|)\to\UU(\g)$ such that 
$$\begin{matrix}\|\g\|&\mapright{\incl}&T(\|\g\|)&{=}&\<\|\g\|\>&\mapright{j_{\|\g\|}}&\TT(|\g|)\\
\lmapdown{\id}&&\lmapdown{p_{\u{\g}}}&&&&\mapdown{q_{\g}}\\
\u{\g}&\lmapright{\incl}&U(\u{\g})&{=}&\<\u{\g}\>&\lmapright{\incl}&\UU(\g)\\
\end{matrix}$$
is commutative.
\end{Theorem}

\msk

\begin{Remark} \label{Lemma- 5.X.} The quotient morphism $q_\g$ respects augmentations
in the sense that $\alpha_\g\circ q_\g = \alpha_{|\g|}$.    
\end{Remark}

\bsk

\section{Appendix: Some facts on weakly complete symmetric Hopf algebras}

\msk

\begin{Definition} \label{primitive}
Let $A$ be a weakly complete symmetric Hopf algebra,
i.e.\ a group object in the monoidal category $(\W,\otimes_W)$ of 
weakly complete vector spaces (see \cite{compbook}, Appendix 7 and
Definition A3.62), with
comultiplication $c\colon A\to A\otimes A$ and coidentity $k\colon A\to \K$. 

(i)  An element $a\in A$ is called
{\it grouplike} if $c(a)=a\otimes a$ and $k(a)=1$. The subgroup
of grouplike elements in the group of units $A^\times$ will be denoted
$\G(A)$. 

\ssk

(ii) An element  $a\in A$ is called {\it primitive}, if
$c(a)= a\otimes 1 + 1\otimes a$. The Lie algebra of primitive elements of 
$\li A$ will be denoted $\P(A)$.                                                              
\end{Definition}

\medskip 

Any weakly complete unital algebra $A$
has an everywhere defined exponential function $\exp\colon \li A\to A^\times$
into the pro-Lie group $A^\times$ of invertible elements defined
as $\exp x=1+x+\frac1{2!}\.x^2+\frac1{3!}\.x^3+\cdots$. As a function
$\exp\colon\li A\to A^\times$ it is the exponential function of the pro-Lie
group $A^\times$ in the sense of pro-Lie groups.

\begin{Theorem} \label{expone} {\rm (Weakly Complete Symmetric Hopf Algebras)}  
If $A$ is a weakly complete symmetric Hopf algebra, then the set
$\G(A)$ of grouplike elements is a closed pro-Lie subgroup of the pro-Lie group
$A^\times$, and the set $\P(A)$ of primitive elements is a closed Lie subalgebra
of the profinite-dimensional Lie algebra $\li A$ and $\exp(\P(A))\subseteq 
\G(A)$ in such a fashion that the restriction and corestriction of $\exp$
is the exponential function $\exp_{\G(A)}:\P(A)\to \G(A)$ of the pro-Lie
group $\G(A)$. 
\end{Theorem}

\nin 
(See e.g.\ \cite{dhtwo}, \cite{compbook}, \cite{probook}.)

\bsk

A simple observation tells us something about the geometry of
the set $\exp A$. Indeed, for
 $t\in\K$ we have $\exp(t\.1+x)=e^t\.\exp x$. Thus

\begin{Remark} \label{trivialtoo} In any weakly complete unital algebra $A$, we have
$$\exp A=\begin{cases}\R_<\.\exp A, \mbox{ if $\K=\R$,}\\
(\C^\times)\.\exp A, \mbox{ if $\K=\C$.}\end{cases}\leqno(*)$$
\end{Remark}

\nin 
Now  we assume that 
{\it $A$ has a coidentity  $\alpha\colon A\to\K$ which is a 
morphism of unital algebras,}
and we call $A$ an {\it augmented algebra}.
We set
$\Ii=\ker \alpha=\{a\in A: \alpha(a)=0\}$.
Then $\Ii$ is a maximal ideal of $A$ and $A/\Ii\cong\K$ and since $\K\.1$
is central we have 
$$A=\K\.1\oplus \Ii\leqno(1)$$
as a direct sum of closed subalgebras.

\begin{Lemma} \label{augmented} Let $A$ be an augmented weakly complete unital algebra.
Then $1+\Ii\subseteq\exp A$.
\end{Lemma}

\begin{Proof} Let $x\in \Ii$ and set $a=1-x\in 1+\Ii$. By  
Lemma \ref{ffactor} (ii)
we find a morphism $\phi\colon \K\<X\>\to A$ of weakly complete unital algebras such that 
$\phi(X)=x$.  By Theorem \ref{monoth} (ii),
in $\K\<X\>$ the element $Y=\log(1-X)=(\log(1-X_f))_{f\in\PP_\K}$ is well defined. Then
$\exp_{\K\<X\>}(Y) =1-X$ and so $a=1-\phi(X)=\phi(1-X)=\phi(\exp_{\K\<X\>}(Y))
=\exp_A(\phi(Y))\in \exp_A(A)$.
\end{Proof}

\begin{Lemma} \label{elementary} If $\K=\R$ then 
$$\R_<\.(1+\Ii)=(\R_<\.1)\oplus \Ii\subseteq\R\.1\oplus\Ii=A,$$
 and if $\K=\C$, then
$$(\C^\times)\.(1+\Ii)=( \C^\times)\,(1\oplus\Ii)=A\setminus\Ii.$$
\end{Lemma}

\begin{Proof} In view of (1) above, the proof is elementary.
\end{Proof}

\begin{Theorem} \label{aug} In any  weakly complete algebra $A$ with 
augmentation $\alpha\colon A\to \K$ let $\Ii\defi\alpha^{-1}(0)$.
\begin{enumerate}[\rm(i)]
\item $\K=\R$: Then $\exp(A)=\alpha^{-1}(\R_<)=(\R_<\.1)\oplus \Ii$.

\item $\K=\C$: Then $\exp A=\alpha^{-1}(\C^\times)=A\setminus\Ii=A^\times$ 
\end{enumerate}
\end{Theorem}

\begin{Proof} By Lemma \ref{augmented} we have $1+\Ii\subseteq\exp(A)$.
By Lemma \ref{trivialtoo} $e^\K\.(1+\Ii)\subseteq\exp A$.
Since $e^\K=\R_<$ for $\K=\R$ and $e^\K=\C^\times$ for $\K=\C$,
Lemma \ref{elementary} completes the proof.
\end{Proof}

\msk

These simple facts complement Theorem \ref{expone}.

\bsk    

\nin{\bf Acknowledgments.} The authors are deeply grateful to the
referee  who has contributed sustantially to the  final form
of this text  in its orthography, typography, and, notably
in the context of Theorems \ref{main-1} and \ref{monoth}, in
 its mathematics.

\ssk 

An essential part of this  text was written while the authors 
were partners in the program {\sc Research  in Pairs} 
at the Mathematisches Forschungsinstitut Oberwolfach MFO 
in the {Black Forest} from February 2 through 22, 2020. The authors
are grateful for  the environment 
and infrastructure of MFO which made this research possible.


\begin{thebibliography}{99}

\bibitem{bour} Bourbaki, N., \flqq Groupes et alg\`ebres de Lie\frqq,
Hermann, Paris, 1971. 

\bibitem{dhtwo} Dahmen, R., and K.H. Hofmann,
{\it The Pro-Lie Group Aspect of Weakly Complete Algebras and
Weakly Complete Group Hopf Algebras},
J. of Lie Theory {\bf29} (2019), 413--455.

\bibitem{dix} Dixmier, J., ``Enveloping Algebras'',
 North-Holland, 1977, xvi+375pp.
 
\bibitem{ham} Hamilton, A., {\it A Poincar\'e-Birkhoff-Witt Theorem
for Profinite Pronilpotent Lie Algebras},
J. of Lie Theory {\bf29} (2019), 611--618.

\bibitem{hilhoflaw} Hilgert, J., K.H. Hofmann, and J.D. Lawson,
{\it Lie Groups, Convex Cones, and Semigroups}, 
Oxford Mathematical Monographs {\bf18}, Clarendon Press Oxford, 1989,
xxxviii+645pp.

\bibitem{hilnee} Hilgert, J., and K.-H. Neeb,
 ``Structure and Geometry of Lie Groups'',
Springer Monographs in Mathematics, Springer NewYork 
etc., 2012, x+744pp.

\bibitem{hofkra} Hofmann, K.H., and L. Kramer, 
{\it On Weakly Complete Group Algebras of Compact Groups},
J. of Lie Theory {\bf30} (2020), 407--424.

\bibitem{hofkraii} Hofmann, K.H., and L. Kramer, 
{\it On Weakly Complete Enveloping Algebras of Profinite-Dimensional
Lie Algebras}, Preprint, Mathematisches Forschungsinstitut Oberwolfach,
 2020, 22pp.

\bibitem{liethree} Hofmann, K.H., and S.A. Morris, 
{\it Sophus Lie's Third Fundamental Theorem and the
Adjoint Functor Theorem,}
J. of Group Theory {\bf8} (2005), 115--133. 

\bibitem{probook}  Hofmann, K.H., and S.A. Morris, ``The Lie Theory of
Connected Pro-Lie Groups,--A Structure Theory for Pro-Lie Algebras,
Pro-Lie Groups, and Connected Locally Compact Groups'', European
Mathematical Society Publishing House, Z\"urich, 2006, xii+663pp.

\bibitem{compbook} Hofmann, K.H., and S.A. Morris, ``The Structure Theory of 
Compact Groups--A Primer for the Student--a Handbook for the Expert'',
4th Edition, De Gruyter Studies in Mathematics {\bf25}, 
Berlin/Boston, 2020, xxvi+1006pp.

\bibitem{hofkraiii} Hofmann, K.H., and S.A. Morris,
{\it Advances in the Theory of Compact Groups and Pro-Lie Groups
in the last Quarter Century}, Axioms {\bf10}, 2021, 190 (13pp.).

\bibitem{hofrup} Hofmann, K.H., and W.A.F. Ruppert,
 {\it Lie Groups and Subsemigroups 
with Surjective Exponential Function}, Memoirs of the Amer.Math.Soc. 
{\bf130}, 1997,
viii+174pp. 

\bibitem{kell} Kelley, J.L., ``General Topology'',
D. van Nostrand, Princeton, New Jersey, 1966.

\bibitem{mac} MacLane, S., ``Categories for the Working Mathematician'',
Graduate Texts {\bf5}, Springer-Verlag, New York etc., 1971, ix+262pp.

\bibitem{montz} Montgomery, D., and L. Zippin, 
 ``Topological Transformation Groups'',
Interscience Publishers, New York, 1955.

\bibitem{most} Mostert, P.S., and A. Shields, {\it Semigroups 
with identity on a manifold},
Trans.  Amer. Math. Soc. {\bf91} (1959), 380-389.

\bibitem{serre-1} Serre, J.-P., 
``Lie Algebras and Lie Groups'', W.A.~Benjamin, 
Inc., New York-Amsterdam, 1965.
   
\bibitem{mathoverflow} Troshkin, M.,
 {\it Group-like elements in a universal
enveloping algebra}, Discussion in mathoverflow, Feb. 4, 2021.
 
\end{thebibliography}
\end{document}